\documentclass{amsart}

\theoremstyle{plain}
\newtheorem{thm}{Theorem}[section]
\newtheorem{lemma}[thm]{Lemma}

\newtheorem{prop}[thm]{Proposition}
\newtheorem{assu}[thm]{Assumption}

\theoremstyle{definition}

\newtheorem{remark}[thm]{Remark}
\newtheorem{remarks}[thm]{Remarks}

\catcode`\@=11
\def\mequal{\mathrel{\mathpalette\@mvereq{\hbox{\sevenrm m}}}} 
\def\@mvereq#1#2{\lower.5\p@\vbox{\baselineskip\z@skip\lineskip1.5\p@
    \ialign{$\m@th#1\hfil##\hfil$\crcr#2\crcr=\crcr}}}
\def\partr#1#2{/\kern-.08333em/_{#1,#2}^{\phantom{.}}}
\def\invpartr#1#2{/\kern-.08333em/_{#1,#2}^{-1}} 
\def\hpartr#1#2{/\kern-.08333em/_{#1,#2}^{h}}
\def\Epartr#1#2{/\kern-.08333em/_{#1,#2}^{E}}

\def\map#1#2#3{{#1}\colon\,{#2}\to{#3}}

\def\newdot{{\kern.8pt\cdot\kern.8pt}}

\def\,{\relax\ifmmode\mskip\thinmuskip\else\thinspace\fi}
\def\{{\relax\ifmmode\lbrace\else $\lbrace$\fi}
\def\}{\relax\ifmmode\rbrace\else $\rbrace$\fi}

\DeclareSymbolFont{rsfs}{U}{rsfs}{m}{n}
\DeclareSymbolFontAlphabet{\mathscr}{rsfs}
\DeclareFontFamily{U}{rsfs}{}
\DeclareFontShape{U}{rsfs}{m}{n}{%
   <5> <6> rsfs5
   <7> rsfs7
   <8> <9> <10> <10.95> <12> <14.4> <17.28> <20.74> <24.88> rsfs10
}{}

\DeclareFontFamily{U}{msb}{}
\DeclareFontShape{U}{msb}{m}{n}{
  <5> <6> <7> <8> <9> gen * msbm
  <10> <10.95> <12> <14.4> <17.28> <20.74> <24.88> msbm10
  }{}
\DeclareSymbolFont{AMSmsb}{U}{msb}{m}{n} 
\SetSymbolFont{AMSmsb}{bold}{U}{msb}{m}{n}
\DeclareSymbolFontAlphabet{\Bbb}{AMSmsb}

\font\sevenrm=cmr7

\newcommand{\C}{{\Bbb C}}
\newcommand{\E}{{\Bbb E}}
\newcommand{\NN}{{\Bbb N}}
\newcommand{\RR}{{\Bbb R}}
\newcommand{\R}{{\Bbb R}}
\newcommand{\OO}{{\Bbb O}}
\newcommand{\PP}{{\Bbb P}}

\newcommand{\SP}{{\mathscr P}}

\def\mathpal#1{\mathop{\mathchoice{\text{\rm #1}}%
   {\text{\rm #1}}{\text{\rm #1}}%
   {\text{\rm #1}}}\nolimits}

\def\Ric{\mathpal{Ric}}
\def\grad{\mathpal{grad}}
\def\Id{\mathpal{Id}}
\def\trace{\mathpal{tr}}

\def\Ito{{\text{\rm It\^o}}}

\def\grad{\mathop{\rm grad}\nolimits}
\def\di{\displaystyle}
\def\f{\frac}
\def\a{\alpha}
\def\b{\beta }
\def\D{\Delta }
\def\d{\delta }
\def\e{\varepsilon }
\def\G{\Gamma }
\def\g{\gamma }
\def\l{\lambda }
\def\n{\nabla }

\def\om{\omega }

\def\s{\sigma }

\def\canu{c^{a}_{\nu}}
\def\CC{{\mathscr C}}

\def\FF{{\mathscr F}}
\def\esp{{\mathbb{E}}}
\def\exsp{{\mathscr{E}}}
\def\H{{\mathbb{H}}}
\def\KK{{\mathscr K}}
\def\lacc{\left\{}
\def\lcr{\left[}
\def\lpa{\left(}
\def\lva{\left|}
\def\Lalb{{\mathscr L}^{\alpha, \beta}}
\def\Lnuq{{\mathscr L}^{\nu, q}}
\def\N{{\mathbb{N}}}
\def\pb{{\mathbb{P}}}
\def\Pmn{{\rm P}^{\mu}_{\nu}}
\def\pnk{\pb^{\nu,k}}
\def\pnak{\pb^{\nu,\alpha, k}}
\def\pbok{\pb^{0,k}}
\def\pbnk{{\bar \pb}^{\nu + 1,k}}
\def\ptnk{{\tilde \pb}^{\nu +1,k}}
\def\psina{\psi^{a, \nu}}
\def\QQ{{\mathbb{Q}}}
\def\racc{\right\}}
\def\rcr{\right]}
\def\rpa{\right)}
\def\rva{\right|}
\def\t\alpha{\tilde{\alpha}}
\def\Tnuq{T^{\nu, q}}

\def\Un{{\bf 1}}
\def\XX{{\mathscr X}}
\def\Xnuq{X^{\nu, q}}
\def\Ynuq{Y^{\nu, q}}

\newcommand{\fin}{\begin{flushright}
                  \mbox{$\Box$}
                  \end{flushright}
                  \noindent}

\begin{document}

\title[Concentration of the Brownian bridge]
      {Concentration of the Brownian bridge on Cartan-Hadamard manifolds with pinched negative sectional curvature}

\author[Marc Arnaudon and Thomas Simon]{Marc Arnaudon and Thomas Simon}

\address{D\'epartement de Math\'ematiques, Universit\'e de Poitiers,T\'el\'eport 2, BP 30179, Boulevard Marie et Pierre Curie, F-86962 Futuroscope-Chasseneuil Cedex, France. {\em E-mail address}: {\tt arnaudon@math.univ-poitiers.fr}}

\address{Equipe d'Analyse et Probabilit\'es, Universit\'e d'\'Evry-Val d'Essonne, Boulevard Fran-\c{c}ois Mitterrand, F-91025 Evry Cedex, France. {\em E-mail address}: {\tt tsimon@univ-evry.fr}}


%
%

\begin{abstract}\noindent  We study the rate of concentration of a Brownian bridge in time one around the corresponding geodesical segment on a Cartan-Hadamard manifold with pinched negative sectional curvature, when the distance between the two extremities tends to infinity. This improves on previous results by A.~Eberle \cite{Eberle:02}, and one of us \cite{Simon:02}. Along the way, we derive a new asymptotic estimate for the logarithmic derivative of the heat kernel on such manifolds, in bounded time and with one space parameter tending to infinity, which can be viewed as a counterpart to Bismut's asymptotic formula in small time \cite{Bismut:84}.
\end{abstract}

\maketitle
\tableofcontents

%
%

\section{Introduction}\label{Intro}
\setcounter{equation}0

Let $M$ be a smooth Cartan-Hadamard manifold with pinched negative sectional curvature, viz. a complete, noncompact, simply-connected $\CC^{\infty}$ Riemannian manifold without boundary, whose all sectional curvatures $\kappa$ satisfy
\begin{equation}
\label{E:1}
-c_2\;\le\; \kappa\;\le\; -c_1
\end{equation}
for some fixed constants $c_2\ge c_1 > 0$. For the sake of concision, later on we will refer to (\ref{E:1}) (resp. to $M$) as the "pinching property"  (resp. as a pinched CH manifold). Let $\rho$ be the Riemannian distance on $M$ and $d \ge 2$ be its dimension. We will make the following further assumption on the curvature tensor of $M$ at infinity:
\begin{assu} 
\label{Christo} For every $\lambda > 0,$ there exists $K_{\lambda} > 0$ such that for every $z\in M$ and every normalized exponential chart centered in $z$ with radius $\lambda$, the Christoffel symbols associated with the Levi-Civita connection in this chart are bounded by $K_{\lambda}$, as well as their derivatives up to order two.
\end{assu}
This assumption holds for example when $M$ is a rank-one symmetric space of the noncompact type (by transitive action of the underlying isometry group) or the universal covering of a compact manifold with pinched negative sectional curvature (by compacity). It seems difficult to find a tractable analogous condition on a global chart diffeomorphic to  $M$ - which exists by Cartan-Hadamard's theorem. For example, some Christoffel symbols associated with Poincar\'e's half-plane model for the hyperbolic plane have an exponential growth. Notice finally that this assumption entails that $\nabla R$ is uniformly bounded on $M$, where $R$ stands for the curvature tensor. However, we got stuck in proving that the converse is true.

For every $x\neq y\in M$, set $\varphi(x,y)=\lacc \varphi(x,y)(t), \; t\in\RR\racc$ for the unit-speed geodesic satisfying $\varphi(x,y)(0)=x$ and $\varphi(x,y)(\rho(x,y))=y$, and 
$$S(x,y)=\lacc\varphi(x,y)(t), \; t\in[0,\rho(x,y)]\racc$$
for the geodesic segment between $x$ and $y$.  Fix $x\in M$ and a unit vector $v\in T_xM$. Define $y(s)=\exp_x(sv)$ for every $s\ge 0$, and consider the $M$-valued Brownian motion  $X(s)=X^{x,y(s)}$ started at $x$ and conditioned to hit the point~$y(s)$ at time~$1$. More precisely, we ask $X(s)$ to solve the It\^o equation 
\begin{equation}
\label{E:4}
d^\n_\Ito X_t(s)\;= \; A(X_t(s))\,dB_t\; + \; V_t(s,X_t(s))\,dt,
\end{equation}
where $B$ is an $\R^m$-valued Brownian motion ($m\ge d$), $A\in\G(\R^m\otimes TM)$ satisfies $A(z)A(z)^{\ast}=\Id_{T_zM}$ for every $z\in M$, $p_t(z,y)$ is the heat kernel on $M$ and  
\begin{equation}
  \label{E:5}
V_t(s,z)\; =\; \grad\log p(1-t,\newdot,y(s))(z).
\end{equation}
Let $y=y(1)$ and $\map{f}{M}{\RR_+}$ be the function $z\mapsto \rho^2(z,\varphi(x,y))$, viz. $f(z)$ is the square of the distance from $z$ to the whole geodesic $\varphi(x,y)$. Consider the process $\lacc Z_t(s)=f(X_t(s)), \; t\ge 0\racc$, and for every $a >0$ the event
$$\Lambda^{x,v,s}_a\; =\; \lacc\sup_{t\in [0,1]} Z_t(s) \; \ge \; a\racc.$$
In the following, we will write $\KK^a_c = (2/c) \log\lpa \cosh c\sqrt{a}\rpa $ for every $a, c >0$. The aim of this paper is to prove the following
\begin{thm}
\label{T1}
Under Assumption \ref{Christo} and with the above notations, for every $a>0$
\begin{eqnarray*}
-\KK^a_{c_2}\;\le\;\liminf_{s\to \infty} s^{-1}\log\PP\lcr \Lambda^{x,v,s}_a\rcr\; \le\; \limsup_{s\to \infty} s^{-1}\log\PP\lcr \Lambda^{x,v,s}_a\rcr\;\le\; -\KK^a_{c_1}, 
\end{eqnarray*}
uniformly in $x\in M$ and $v\in T_xM$ such that  $\|v\|=1$. Besides, the same result holds in replacing $\varphi(x,y)$ by $S(x, y)$ in the definition of $\Lambda^{x,v,s}$.
\end{thm}
This result means that the Brownian bridge "concentrates" around the geodesic line - resp. the geodesic segment - joining its two extremities when the distance between the latter tends to infinity, and extends to pinched CH manifolds the main theorem of~\cite{Simon:02}, which established the result on the real hyperbolic plane with constant sectional curvature -1 (in this case one finds then an exact limit given by $\KK^a_1$). Recall that originally, a weak version of this concentration result had been obtained by Eberle \cite{Eberle:02}, providing the key-step in the construction of a counterexample for the existence of a spectral gap on the loop space over a compact Riemannian manifold. 

The main argument used in this paper to obtain Theorem \ref{T1} is entirely different from the techniques developed in \cite{Eberle:02} and \cite{Simon:02}, where the Brownian bridge was rather considered as an $h$-transform of Brownian motion. Here, we choose to work directly on the SDE (\ref{E:4}), and the main point consists in obtaining the following limit theorem for its drift coefficient when $s\to \infty$:
\begin{equation}
\label{limit}
\lim_{s\to +\infty}s^{-1}V_t(s,X_t(s))\; = \;\dot\varphi(X_t(s),y(\infty))(0).
\end{equation}
Indeed, once (\ref{limit}) is obtained, a simple application of It\^o's formula to the process $Z$, combined with Alexandrov-Toponogov's triangle comparison theorem  and the comparison theorem for real SDE's shows that when $s\to\infty$, a.s. $Z_t$ lies roughly between the solutions of the SDE's
$$Y^i_t\; =\; 2\int_0^t \sqrt{Y^i_u} dB_u\; -\; 2s\int_0^t \sqrt{Y^i_u} \tanh \lpa c_i \sqrt{Y^i_u}\rpa du\; +\; k_i t, \quad i = 1, 2$$
where $k_1, k_2$ are positive constants. An asymptotic analysis of these latter diffusions of the Cox-Ingersoll-Ross type, performed with the help of stochastic calculus and first passage time techniques, delivers then the required lower and upper exponential speeds of convergence $-s\KK^a_{c_i}$, $i = 1, 2$. By the way, we remark that these speeds of convergence can be computed in integrating from $0$ to $a$ the functions $2b_i(x)/a^2(x)$, where $a(x) = 2\sqrt{x}$ is the diffusion coefficient and $b_i(x) = -2s\sqrt{x}\tanh\lpa c_i \sqrt{x}\rpa$ is the dominating drift coefficient of the corresponding diffusion. However, we could not find any sensible geometrical explanation of this computation.

Up to technical details - which are a bit reminiscent to those of Eberle's paper, the limit theorem (\ref{limit}) is actually a direct consequence of the following logarithmic derivative estimate of the heat kernel on $M$, when $s\to +\infty$:
\begin{equation}
\label{limit2}
s^{-1}\grad\log p_t(\newdot,y)(x)\; \to \; t^{-1}\dot\varphi(x,y)(0).
\end{equation}
This estimate, which may be interesting by itself as a pendant to Bismut's celebrated estimate for $\grad\log p_t(\newdot,y)(x)$ in small time - see Theorem 3.8 in \cite{Bismut:84}, is rather easy to obtain analytically on real hyperbolic spaces or rank-one noncompact symmetric spaces, because of the existence of (more or less) closed formulae thereon. The situation is however much more complicated on general pinched CH manifolds. To achieve our proof, we use then probabilistic arguments relying on a "filtered" integration by part formula for the heat kernel \cite{Thalmaier:97} \cite{Thalmaier-Wang:98}, and a suitable development in local coordinates where we can perform large deviation estimates and apply Varadhan's lemma. In the end, we do obtain (\ref{limit2}) in full generality on $M$, but unfortunately we need the Assumption \ref{Christo} to obtain the {\em uniform} convergence in $x$, which is crucial to get (\ref{limit}). This explains the restriction on $M$ in the statement of Theorem \ref{T1}. At the end of the paper, we provide an example where (\ref{limit2}) may fail in the absence of Assumption \ref{Christo}.

In addition to being more general, we feel that our proof is more transparent than the ones in \cite{Eberle:02} and \cite{Simon:02}, even though it would be quite interesting to see if Martin boundary techniques could also apply on general CH manifolds to study this concentration phenomenon. Yet another approach could be the following: observing by Riemannian comparison that if $M_1$, $M_2$ are two CH manifolds such that $\sup_{x\in M_2} \kappa_2(x)\; \le \; \inf_{x\in M_1} \kappa_1(x)$ - with the above notations, then for any given geodesics $\gamma_i\subset M_i$ it is possible to construct two Brownian motions $X^i$ starting from $x_i\in\gamma_i$, $i =1,2$, such that a.s. $\rho(X^2_t, \gamma_2)\ge  \rho(X^1_t, \gamma_2),$ one may wonder if conditioning both $X^i$'s to go back to $\gamma_i$ in time 1 should not force $X^2$ to stay closer to $\gamma_2$ than $X^1$ to $\gamma_1$ in the meantime. Roughly, this would then prove Theorem \ref{T1} provided the result is already known on real hyperbolic spaces, because of (\ref{E:1}). Unfortunately, we could not give a rigorous approach to these simple heuristics relying only on the constant curvature case and Riemannian comparison theorems.

\section{The case of real hyperbolic spaces}\label{Hyperbolic}\setcounter{equation}0

In this section we generalize the main result of \cite{Simon:02} to all real hyperbolic spaces $\H^d_c(\R)$, $d\ge 2$, with constant sectional curvature $-c < 0$. As in \cite{Simon:02}, the first step consists in estimating the deviations from the origin of a family of diffusions of the Ornstein-Uhlenbeck type with big negative drift. Then we prove the logarithmic derivative estimate of the heat kernel and two further estimates, which entail together with It\^o's formula that the process $\lacc Z_t, \; t\geq 0\racc$ becomes very close to these diffusions when $s\to +\infty$. The proof of Theorem \ref{T1} follows then simply from the comparison theorem for one-dimensional stochastic differential equations.

\subsection{Asymptotics of first-passage times for CIR-type diffusions}

We begin with an extension of the Proposition of \cite{Simon:02}, showing that the limit constant therein actually does not depend of the dimension. For every $\nu\in\R$ and $c, k > 0$, let $Y^{\nu, c, k}$ be the solution to the SDE
\begin{equation} \label{SDE}
Y^{\nu, c, k}_t\; =\; 2\int_0^t \sqrt{Y^{\nu, c, k}_s} dB_s\; -\; 2\nu\int_0^t \sqrt{Y^{\nu, c, k}_s} \tanh \lpa c \sqrt{Y^{\nu, c, k}_s}\rpa ds\; +\; kt
\end{equation}
where $\lacc B_t, \; t\ge 0\racc$ is a standard linear Brownian motion. By analogy with Bessel diffusions, we see that $(\ref{SDE})$ has a unique strong solution which is positive for every $t > 0$. In the following we will set $\pb^{\nu, c, k}$ for the law of $Y^{\nu, c, k}$. If $\lacc X_t, \; t\ge 0\racc$ is the canonical process, let $\lacc \FF_t, \; t\ge 0\racc$ be the canonical completed filtration, and $T_a$ be the first hitting time of $X$ at level $a > 0$: $T_a \; = \; \inf\lacc t > 0 \; / \; X_t = a \racc.$ Notice that under $\pb^{0, c, k}$, $X$ is the square of a Bessel process of dimension $k$, which we will be sometimes denote by $X^k$ when no confusion is possible. 

\begin{prop} \label{1dim1} For every $a, c, k, t > 0$, 
\begin{eqnarray*}
\lim_{\nu\uparrow +\infty} \nu^{-1}\, \log \pb^{\nu, c, k}\lcr  T_a < t \rcr & = & - \KK^a_c
\end{eqnarray*} 
\end{prop}

\noindent
\begin{proof} We first notice that it suffices to consider the case $c = 1$, the general case $c > 0$ following from a straightforward scaling argument in considering the process $t\mapsto c^2 Y^{\nu, c, k}_{t/c^2}$. Setting $\pnk = \pb^{\nu, 1, k}$ for simplicity, we will show that 
$$- \KK^a_1\;\leq\;\liminf_{\nu\uparrow +\infty} \nu^{-1}\, \log \pnk\lcr  T_a < t \rcr\; \leq \;\limsup_{\nu\uparrow +\infty} \nu^{-1}\, \log \pnk\lcr  T_a < t \rcr\; \leq \; - \KK^a_1.$$

\noindent
{\bf Proof of the lower limit}. By Girsanov's theorem and the fact that $\lacc T_a < t\racc\in \FF_{T_a}$, we can write
\begin{eqnarray*}
\pnk\lcr T_a < t\rcr & = & \pbok\lcr T_a < t; \; L^{\nu,k}_{T_a}\rcr
\end{eqnarray*}
where
$$L^{\nu,k}_{T_a} \; = \; \exp \lcr -\nu \int_0^{T_a}\frac{\tanh\sqrt{X_s}}{2\sqrt{X_s}} \,(dX_s - kds) -\frac{\nu^2}{2}\int_0^{T_a}\tanh^2\sqrt{X_s} \,ds\rcr.$$
On the other hand, It\^o's formula yields
\begin{eqnarray*}
 -\nu \int_0^{T_a}\frac{\tanh\sqrt{X_s}}{2\sqrt{X_s}} \, dX_s & = & -(2\nu +1) \log\,\cosh \sqrt{a}  \; + \; (\nu +1) \int_0^{T_a}\frac{\tanh\sqrt{X_s}}{2\sqrt{X_s}} \, dX_s\\
 & & +(\nu + 1/2) \int_0^{T_a}\lpa 1 - \frac{\tanh\sqrt{X_s}}{\sqrt{X_s}}  - \tanh^2\sqrt{X_s}\rpa \, ds
\end{eqnarray*}
and we get, after some rearrangements,
\begin{eqnarray}
\label{Girsanov}
\pnk\lcr T_a < t\rcr & = & \lpa \cosh \sqrt{a}\rpa^{-(2\nu +1)}\pbnk\lcr T_a < t; \; M^{\nu,k}_{T_a}\rcr,
\end{eqnarray}
where we set ${\bar \pb}^{\nu,k} = \pb^{-\nu,k}$ and 
\begin{eqnarray*}
M^{\nu,k}_{T_a} & = & \exp\lcr(\nu + 1/2) \int_0^{T_a}\lpa 1 + (k-1)\frac{\tanh\sqrt{X_s}}{\sqrt{X_s}}\rpa \, ds\rcr.
\end{eqnarray*}
We will now prove that $\pbnk\lcr T_a < t\rcr $ tends to 1 as $\nu \to +\infty$, which is sufficient to obtain the lower limit, because $M^{\nu,k}_{T_a} \ge 1$ a.s. Setting $c_a = \tanh \sqrt{a}/\sqrt{a}$, we see by comparison that if $\ptnk$ stands for the law of the solution to the SDE
$$X_t\; =\; 2\int_0^t \sqrt{X_s} dB_s\; + \;2c_a(\nu+1)\; \int_0^t X_s \, ds\; +\; kt,$$
then $\pbnk\lcr T_a < t\rcr \ge \ptnk\lcr T_a < t\rcr$ for every $t >0$. Under $\ptnk$, we recognize in $X$ the well-known Cox-Ingersoll-Ross (CIR) process, which can be reconstructed from the square Bessel process $X^k$ by deterministic time change:
\begin{equation}
\label{timechange}
\lacc X_t, \; t\ge 0\racc\; \stackrel{d}{=}\; \lacc e^{2\canu t} X^k_{\psina_t},\; t\ge 0\racc,
\end{equation}
where we set $\canu = (\nu+1)c_a$ and $\psina_t = \lpa 1 - e^{-2\canu t} \rpa/2\canu
$. Using the scaling property of $X^k$, this entails 
\begin{eqnarray}
\label{bessel}
\ptnk\lcr T_a < t\rcr & \ge & \pb\lcr X^k_{\psina_t}> a e^{-2\canu t} \rcr\nonumber\\
& \ge & \pb\lcr X^k_1 >  2\canu ae^{-2\canu t}/(1 - e^{-2\canu t})\rcr,
\end{eqnarray}
which completes the proof of the lower limit, because $X^k_1$ does not weight $\{0\}$ for every $k >0$, and since $2\canu ae^{-2\canu t}/(1 - e^{-2\canu t})\to 0$ when $\nu \to +\infty$.\\

\noindent
{\bf Proof of the upper limit}. We will use a different method, relying on Feller's spectral theory \cite{Feller:54}. Actually, the case $k=1$ was already proved in the same way in the Proposition of \cite{Simon:02}, with the help of Legendre functions. But the situation is a bit different when $k \neq1$, because the underlying second order differential equation has then a new singularity at zero - see however the following Remark 2.2.(c) for the case $k=3$. Fix $a, t > 0$, $x\in ]0, a[$, and set $\pnk_x$ for the law of the solution to (\ref{SDE}) starting from $x$, whence $\pnk\lcr T_a < t\rcr \le \pnk_x\lcr T_a < t\rcr$ by comparison. Assuming without loss of generality that $k\ge 2$, setting $q = (k-1)/2$ and $\Xnuq$ for the unique (positive) solution to the SDE 
\begin{eqnarray*}
\Xnuq_s& = & \sqrt{x} \; +\; B_s\; - \;\nu \int_0^s\!\! \tanh \Xnuq_u\, du \; +\; q \int_0^s\!\! \frac{1}{\Xnuq_u}\,du,
\end{eqnarray*}
we see from It\^o's formula that under $\pnk_x$ we have a.s. $T_a = \inf \lacc s > 0, \; \Xnuq_s= \sqrt{a}\racc$. Hence, another comparison argument yields 
\begin{equation}\label{tmuaq}
\pnk\lcr T_a < t\rcr \;\le \;\pnk_x\lcr T_a < t\rcr \;\le\; \pb\lcr \Tnuq_a < t\rcr,
\end{equation}
 where $\Tnuq_a = \inf \lacc s > 0, \; \Ynuq_s = \sqrt{a}\racc$ and $\Ynuq$ solves
\begin{eqnarray*}
\Ynuq_s & = & \sqrt{x} \; +\; B_s\; - \;\nu \int_0^s\!\! \tanh \Ynuq_u\, du \; +\; q \int_0^s\!\! \coth \Ynuq_u\, du.
\end{eqnarray*}
Notice now that since $k \ge 2$ and according to Feller's classification - for which we use Mandl's terminology, see \cite{Mandl:68} pp. 13, 24-25 and 67 - the diffusion $\Ynuq$ has a natural boundary at $+\infty$ and an entrance boundary at 0. Hence, according to Feller's spectral theory - see Theorem 4 p. 11 in \cite{Feller:54} or Lemma 3 p. 62 in \cite{Mandl:68}, for every $\lambda > 0$ the Laplace transformation $\esp [ e^{-\lambda \Tnuq_a} ]$ is given by the value at $z = \sqrt{x}$ of the unique solution over $]0, \sqrt{a}]$ to the differential equation 
\begin{eqnarray} \label{ODE}
\Lnuq f (z) \; - \; 2\lambda f(z) & = & 0
\end{eqnarray}
satisfying $f\lpa\sqrt{a}\rpa = 1$ and such that $f(0+)$ and $\Lnuq f(0+)$ exist, where we set
\begin{eqnarray*}
\Lnuq f (z) & = & f''(z)\; - \; 2\nu \tanh z f'(z) \; +\; 2q \; \coth z f'(z)
\end{eqnarray*}
for every $z > 0$ and smooth functions $f$. The above operator $\Lnuq$ is in fact well-known from harmonic analysis as a Jacobi operator - see Section 2 in \cite{Koornwinder:84}. Setting $\mu = \sqrt{(\nu-q)^2 + 2\lambda}$, we see from the boundary conditions and Formula (2.7) in \cite{Koornwinder:84} that $\esp [e^{-\lambda T^{\nu, q}_a} ]$ equals
$$\frac{\lpa \cosh \lpa \sqrt{x}\rpa\rpa^{\nu + \mu -q} {\rm F} \lpa  (q - \nu - \mu)/2 \; , \; (q + 1 + \nu -\mu)/2\; ; \; q + 1/2 \; ; \; \tanh^2\lpa \sqrt{x}\rpa \rpa}{\lpa \cosh \lpa \sqrt{a}\rpa\rpa^{\nu + \mu -q}  {\rm F} \lpa  (q - \nu - \mu)/2 \; , \; (q + 1 + \nu -\mu)/2\; ; \; q + 1/2 \; ; \; \tanh^2\lpa \sqrt{a}\rpa \rpa},$$
where ${\rm F}$ stands for Gauss' hypergeometric function. By an asymptotic expansion of the latter when its first parameter is (negatively) large - see e.g. \cite{Magnus:66} p. 56, this entails finally
$$\limsup_{\nu\to +\infty} \nu^{-1} \log \pb \lcr  T^{\nu, q}_a < t \rcr\; \leq\; \lim_{\nu\to +\infty} \nu^{-1} \log \esp \lcr e^{-\lambda T^{\nu, q}_a} \rcr\; =\; \KK^x_1 - \KK^a_1,$$
where the first inequality is an immediate consequence of the Markov inequality. Using (\ref{tmuaq}) and letting $x$ tend to 0 completes now the proof of the upper limit. 

\end{proof}

\begin{remarks} (a) From the proof of the lower limit, we notice that the exponential speed of convergence $-\nu \KK^a_c$ emerges naturally in integrating from $0$ to $a$ the function $2b(x)/a^2(x)$, where $a(x) = 2\sqrt{x}$ is the diffusion coefficient and $b(x) = -2\nu\sqrt{x}\tanh\lpa c \sqrt{x}\rpa$ is the dominating drift coefficient. Actually, the statement of Proposition \ref{1dim1} probably holds for a more general class of diffusions of the square Ornstein-Uhlenbeck type, with big negative drift. However, we notice that the change of measure given by (\ref{Girsanov}) seems useless to obtain the upper limit. Indeed, for example when $k = 1$, it is possible to compute
$$\bar{\esp}^{\nu + 1, 1} \lcr M_{T_a}^{\nu, 1}\rcr\; = \; \lpa \cosh \sqrt{a}\rpa^{2\nu},$$ 
so that we strongly need to consider the event $\lacc T_a < t\racc$ in the analysis of the upper limit. If we could prove {\em a priori} that $\lim_{\nu\uparrow +\infty} \nu^{-1}\, \log \pb^{\nu, c, k}\lcr  T_a < t \rcr$ exists and does not depend on $t$, then this would give a quicker proof of the upper bound without special functions, thanks to the immediate inequality
$$\log\pbnk\lcr T_a < t; \; M^{\nu,k}_{T_a}\rcr\; \leq \; k(\nu + 1/2)t.$$
Unfortunately, we could not find any ergodic theoretical argument removing the dependence on $t$ at the limit, under the logarithmic scale.

\vspace{2mm}

\noindent
(b) In the case $q =0$ (i.e. $k = 1$), the proposition was already proved in \cite{Simon:02} with the help of Legendre functions. Taking this for granted, one can give the following alternative proof of the proposition in the case $k > 1$ - which is the only situation relevant to our further purposes. With the above notations,  it is sufficient to prove that
\begin{equation}
\label{comparison}
 X_s^{\nu,q}\;\le\; X_s^{\nu',0}+\a_0,
\end{equation}
where $\di \nu'=\nu- q/(\a_0\tanh \a_0)$. Indeed, clearly the converse inequality $X_s^{\nu,0}\le X_s^{\nu,q}$ holds a.s. and then we can reason exactly as above. But by It\^o-Tanaka's formula, we have 
$$X_s^{\nu,0}\; =\; B_s-\nu\int_0^s\tanh X_u^{\nu,0}\,du+L_s^{\nu,0},$$
where $L_s^{\nu,0}$ is the local time at $0$ of $X_s^{\nu,0}/2$. On the other hand, the process $X_s^{\nu',\a_0,0}=X_s^{\nu',0}+\a_0$ has for drift 
$$-\nu'\tanh\left(X_s^{\nu',\a_0,0}-\a_0\right)\,ds+dL_s^{\nu,0}$$
whereas when  $X_s^{\nu,q}>0$, the process $X_s^{\nu,q}$ has for drift 
$$\left(-\nu\tanh X_s^{\nu,q}+\f{k-1}{2X_s^{\nu,q}}\right)\,ds.$$
Since for positive time $X_s^{\nu',\a_0, 0}$ is always larger than or
equal to $\a_0$, and since the local time $L_s^{\nu,0}$ is nondecreasing, it
is sufficient to prove that  for every $x\ge  \a_0$,
$$-\nu'\tanh\left(x-\a_0\right)> -\nu\tanh x+(k-1)/2x.$$
A sufficient condition is clearly $\di(\nu-\nu')\tanh x\ge (k-1)/2x.$
Observing that the left hand side is an increasing function of $x$ and
the right hand side is a decrasing function of $x$, a sufficient
condition becomes $\di (\nu-\nu')\tanh \a_0\ge (k-1)/2\a_0,$ so that finally, letting 
$$\nu'\;=\;\nu-\f{k-1}{2\a_0\tanh \a_0},$$
we obtain \eqref{comparison}.

\vspace{2mm}

\noindent
(c) In the case $q =1$ (i.e. $k = 3$), the equation (\ref{ODE}) can be solved in a different way. Using the substitution $f(z) = (\cosh z)^{\nu} (\sinh z)^{-1}g(\tanh z)$ where $g : (-1, 1) \to \R$ is some unknown function, yields namely the following equation for $g$:
$$(1 - z^2)\, g''(z) \; - \; 2z\, g'(z) \; + \; \lpa \nu(\nu +1) - \frac{(\nu-1)^2
+ 2\lambda}{1- z^2}\rpa g(z) \; = \; 0,$$
which is Legendre's differential equation on the cut. Hence, the general solution to (\ref{ODE}) has the form 
$$f(z) \; = \; {(\cosh z)}^{\nu}(\sinh z)^{-1}\lcr A\,\Pmn (\tanh z)\; + \; B\,\Pmn
(-\tanh z) \rcr$$
for two unknown constants $A$ and $B$, where $\mu = \sqrt{(\nu-1)^2 + 2\lambda}$
and $\Pmn$ stands for the Legendre function of the first kind. From the boundary conditions and the third formula p. 167 in \cite{Magnus:66}, we deduce that $\esp [e^{-\lambda T^{\nu, 1}_a} ]$ equals
$$\frac{\lpa \cosh \lpa \sqrt{x}\rpa\rpa^{\nu + \mu -1} {\rm F} \lpa  1/2 - (\nu + \mu)/2 \; , \; 1 + (\nu -\mu)/2\; ; \; 3/2 \; ; \; \tanh^2\lpa \sqrt{x}\rpa \rpa}{\lpa \cosh \lpa \sqrt{a}\rpa\rpa^{\nu + \mu -1}  {\rm F} \lpa  1/2 - (\nu + \mu)/2 \; , \; 1 + (\nu -\mu)/2\; ; \; 3/2 \; ; \; \tanh^2\lpa \sqrt{a}\rpa \rpa},$$
which is of course the same formula as above, for $q=1$. Since $q + 1/2\not\in\lacc 1/2, 3/2\racc$ when $q \not\in\lacc 0, 1\racc$, and recalling the formulae p. 167 in \cite{Magnus:66}, it seems that apart from the regular case $k=1$, the resolution of (\ref{ODE}) with Legendre functions is only possible when $k =3$. We could not find a sensible explanation of this fact.

\vspace{2mm}

\noindent
(d) The above operator $\Lnuq$ plays a central r\^ole in harmonic analysis on rank-one non-compact symmetric spaces, since for suitable choices of $\nu, q, \lambda$ the odd solutions to (\ref{ODE}) yield all the spherical functions on such spaces - see Part 4 in \cite{Koornwinder:84} for much more on this topic. However, the connection between these spherical functions and our equation (\ref{SDE}) is only apparent. Namely, Jacobi operators related to spherical functions on rank-one groups have the form
$$\Lalb f(z)\; =\; f''(z) \; +\; ((2\alpha +1) \coth z \, +\, (2\beta +1) \tanh z) f'(z)$$
with $\alpha > \beta\ge 0$ - see (3.4) in \cite{Koornwinder:84} or our table in  Section 3.1. below, whereas in (\ref{SDE}) our coefficient $q$ before $\coth$ can be neglected in the analysis.
\end{remarks}

In the following, it will be important to consider the perturbation of the above SDE (\ref{SDE}) by some parameter $\alpha \in\R$:
\begin{eqnarray}\label{SDE2}
Y^{\nu, \alpha, c, k}_t & = & 2\int_0^t \sqrt{Y^{\nu, \alpha, c, k}_s} dB_s \; + \; kt\\
& - &  2\nu\int_0^t \sqrt{Y^{\nu, \alpha, c, k}_s} \lpa \tanh \lpa c \sqrt{Y^{\nu, \alpha, c, k}_s} \rpa+ \alpha \rpa ds.\nonumber
\end{eqnarray}
Again, this equation has a unique strong solution which is positive for every $t > 0$. We will set $\pb^{\nu, \alpha, c, k}$ for the law of $Y^{\nu, \alpha, c, k}$ and use the same notations as above for $T_a$ and the canonical process. The proof of the following proposition is very similar to the one above, but requires heavier notations and so we wrote it down separately, for the sake of clarity.

\begin{prop} \label{1dim2} For every $a, c, k, t > 0$, 
$$\lim_{\alpha \to 0}\lpa \liminf_{\nu\uparrow +\infty} \nu^{-1}\log \pb^{\nu, \alpha, c, k}\lcr  T_a < t \rcr \rpa\; = \;\lim_{\alpha \to 0}\lpa\limsup_{\nu\uparrow +\infty} \nu^{-1} \log \pb^{\nu, \alpha, c, k}\lcr  T_a < t \rcr\rpa$$
and the common limit equals $- \KK^a_c $.
\end{prop}

\noindent
\begin{proof} By the same scaling argument as above, it suffices to consider the case $c = 1$ and we will set $\pnak = \pb^{\nu, \alpha, 1, k}$. We begin with the lower limit:
$$\liminf_{\alpha \to 0}\lpa \liminf_{\nu\uparrow +\infty} \nu^{-1}\log \pb^{\nu, \alpha, k}\lcr  T_a < t \rcr \rpa\; \geq \; - \KK^a_1,$$
and we notice that thanks to Proposition \ref{1dim1} and a comparison argument, it suffices to consider the situation where $\alpha > 0$ and $\alpha \downarrow 0$.  
Suppose first that $k > 1$. Applying It\^o's formula to $\log \cosh \sqrt{X_{T_a}} + \alpha \sqrt{X_{T_a}}$ and reasoning exactly as in Proposition \ref{1dim1} yields
$$\pnak\lcr T_a < t\rcr\; =\;\lpa e^{\alpha \sqrt{a}} \cosh \sqrt{a}\rpa^{-(2\nu + 1)}{\bar \pb}^{\nu+1,\alpha,k} \lcr T_a < t; \; M^{\nu,\alpha, k}_{T_a} \rcr,$$
where similarly we set 
${\bar \pb}^{\nu,\alpha, k} = \pb^{-\nu, \alpha, k}$ and $M^{\nu,\alpha, k}_{T_a}$ is given by
$$\exp\lcr(\nu + 1/2) \int_0^{T_a}\lpa 1 +  2\alpha\tanh\sqrt{X_s} + \alpha^2 + (k-1)\frac{\lpa \tanh\sqrt{X_s} +\alpha\rpa}{\sqrt{X_s}}\rpa \, ds\rcr.$$
Since $k \ge 1$, we have again a.s. $M^{\nu,\alpha, k}_{T_a} \geq 1$. Besides, the comparison 
$${\bar \pb}^{\nu+1,\alpha,k} \lcr T_a < t\rcr\;\ge\;  \ptnk\lcr T_a < t\rcr$$
still holds because $\alpha > 0$, so that we can finish the proof of the lower limit exactly as in Proposition \ref{1dim1}.

When $k < 1$ - this case is actually irrelevant to our further purposes but we treat it for completeness, the above method fails because the Wiener integral $\int_0^s X_s^{-1/2}\, ds$ diverges, as for Bessel diffusions. However, we can reinterpret (\ref{SDE2}) as (\ref{SDE}) driven by some drifted Brownian motion $\lacc B^{\nu, \alpha}_s = B_s - (\alpha\nu) s, \; s\ge 0\racc$, and the Cameron-Martin formula yields 
\begin{eqnarray*}
\pnak \lcr T_a < t\rcr & = & \pnk\lcr T_a < t; \; e^{-\alpha\nu B_t - \alpha^2\nu^2t/2}\rcr.
\end{eqnarray*}
Introduce now $K > 0$ and suppose $\nu > K/t$. For $\alpha$ small enough, we first get
\begin{eqnarray*}
\pnak \lcr T_a < t\rcr & \ge & \pnak \lcr T_a < K/\nu\rcr\\ 
& \ge & e^{-2\nu\alpha K}\pnk\lcr T_a < K/\nu ; \; B_{K/\nu} \le K\rcr\\
& \ge & e^{-2\nu\alpha K}\lpa\pnk\lcr T_a < K/\nu\rcr - \pb\lcr B_{K/\nu} > K\rcr\rpa\\
& \ge & e^{-2\nu\alpha K} \lpa \lpa \cosh \sqrt{a}\rpa^{-(2\nu +1)}\pbnk\lcr T_a < K/\nu \rcr- {\rm Erfc} \lpa \sqrt{K \nu}\rpa\rpa.
\end{eqnarray*}
where ${\rm Erfc}$ stands for the Gaussian error function \cite{Magnus:66}. Besides, using (\ref{bessel}) and choosing $K$ big enough, we obtain  
\begin{eqnarray*}
\pbnk\lcr T_a < K/\nu \rcr & \ge & \pb\lcr X^k_1 >  \lpa 4c_a ae^{-2Kc_a}\rpa\nu\rcr.
\end{eqnarray*}
Recalling now that the density of $X^k_1$ over $\R^+$ is given by the function 
$$x\mapsto \lpa 2\Gamma (k/2)\rpa^{-1}(x/2)^{k/2-1} e^{-x/2}$$
see e.g. Corollary XI.1.4 in \cite{Revuz-Yor:99}, and plugging the two above inequalities together, we see that when $K$ is big enough, then 
$$\liminf_{\nu\to +\infty} \nu^{-1}\log \pb^{\nu, \alpha, k}\lcr  T_a < t \rcr \; \geq \; - \KK(a,\alpha, K),$$
where $\KK(a, \alpha, K)$ is some constant tending to $\KK^a_1$ when $\alpha \downarrow 0$ and then $K \uparrow +\infty$. This completes the proof of the lower limit for $k \le 1$. To prove the upper limit:
$$\limsup_{\alpha \to 0}\lpa \limsup_{\nu\uparrow +\infty} \nu^{-1}\log \pb^{\nu, \alpha, k}\lcr  T_a < t \rcr \rpa\; \leq \; - \KK^a_1,$$
we will use another comparison argument. First we can assume $k\geq 2$ without loss of generality and again, we only need to consider the case where $\alpha < 0$ and $\alpha \uparrow 0$. Set $\t\alpha = - \arg\tanh \alpha > 0$ with $\alpha$ small enough, fix $x\in]-\t\alpha^2, a[$ and let $Y^{\nu, \alpha, q}$ be the unique positive solution to the SDE 
$$Y^{\nu, \alpha, q}_s\; =\; \sqrt{x}\; +\; B_s\; -\; \nu\int_0^s \lpa \tanh \lpa  Y^{\nu, \alpha, q}_u\rpa -\tanh \t\alpha \rpa du\; +\; q\int_0^s \frac{du}{Y^{\nu, \alpha, q}_u},$$
with $q = (k-1)/2 \geq 1$. Recall that $\pb^{\nu, \alpha, k}\lcr  T_a < t \rcr \leq  \pb\lcr  T^{\nu, \alpha, q}_a < t \rcr$, where we set $T^{\nu, \alpha, q}_a = \inf\lacc s >0, \; Y^{\nu, \alpha, q}_s = \sqrt{a}\racc$. Using the inequality $(\tanh x - \tanh \t\alpha) \ge c_{\alpha} \tanh (x -\t\alpha),$ which holds uniformly on $x\in \lcr \t\alpha, \sqrt{a} \rcr$ for some constant $c_{\alpha} < 1$ tending to 1 as $\alpha \uparrow 0$, we can compare $Y^{\nu, \alpha, q}$ with $Z^{\nu, \alpha, q}$ solution of
$$Z^{\nu, \alpha, q}_s\; =\; \sqrt{x}\; + \; B_s\; -\; c_{\alpha}\nu\int_0^s \lpa \tanh \lpa Z^{\nu, \alpha, q}_u - \t\alpha\rpa\rpa\,du\; +\; q\int_0^s \frac{du}{\lpa Z^{\nu, \alpha, q}_u -\t\alpha\rpa}$$
(which remains a.s. above the level $\t\alpha >0,$ because $q\ge 1$), and we obtain $\pb\lcr S^{\nu, \alpha, q} < t\rcr \le\pb\lcr S^{\nu, \alpha, q} < t\rcr$ with the notation $S^{\nu, \alpha, q}_a = \inf\lacc s >0, \; Z^{\nu, \alpha, q}_s = \sqrt{a}\racc$. Introducing the process ${\tilde Z}^{\nu, \alpha, q}_s = Z^{\nu, \alpha, q}_s - \t\alpha$ for every $s >0$ and setting $\tilde{x} = (\sqrt{x} - \t\alpha)^2$ and $\tilde{a}= (\sqrt{a} - \t\alpha)^2$, we finally get
$$\limsup_{\nu\uparrow +\infty} \nu^{-1} \log \pb^{\nu, \alpha, k}\lcr  T_a < t \rcr \; \le\; \KK^{\tilde{x}}_1 - \KK^{\tilde{a}}_1,$$  
which finishes the proof of the upper limit in letting $\alpha$, and then $x$, tend to 0. 

\end{proof}     

\subsection{Three further estimates}

In this subsection we establish three crucial estimates which will allow us later on to reduce the original problem to the above asymptotic study for CIR-type processess. The first estimate is fairly straightforward:

\begin{lemma} \label{Laplace} Let $\phi$ be a geodesic line in $\H^d_c$. Setting $g(z) = \rho (z, \phi)$ and $f(z) = g(z)^2$ for every $z\in\H^d_c$, then the following inequalities hold 
$$2 \; \leq\; \Delta f(z)\; \leq\; d \; + \; g(z)$$ 
uniformly in $\H^d_c$.
\end{lemma}

\noindent
\begin{proof}  We first notice that by definition $\lva\lva\grad g\rva\rva = 1$, whence
$$\Delta f(z)\; = \; 2 \;+ \; g(z) \Delta g(z).$$
It remains to estimate $\Delta g(z)$, and this is done in choosing for $\H^d_c\subset\R^d$ a half-space model $\lacc z_d > 0 \racc$, and for $\phi$ the $z_d$-axis. For $z = (z_1, \ldots, z_d)$ we then have $g(z) = c^{-1}\arg \sinh (r_{d-1}/z_d)$, where we set $r_{p} = (z_1^2 + \ldots + z_{p}^2)^{1/2}$ for $p = 1 \ldots d $. A direct computation yields 
$$\Delta g(z) \; = \; c z_d^2\lpa \partial^2_{z_1} + \ldots + \partial^2_{z_d}\rpa\lpa c^{-1}\arg \sinh (r_{d-1}/z_d)\rpa \; = \; \frac{(d-2)z_d^2 + r_{d-1}^2}{r_{d-1}r_d},$$
which entails
$$
0\; \le\; g(z)\Delta g(z)\;\le \; g(z)\lpa \frac{(d-2)z_d}{r_{d-1}} + 1\rpa\; \le \; (d-2) \; + \; g(z)$$
and completes the proof.

\end{proof}

\begin{remarks} (a) If $\phi$ is a geodesical segment, one can prove that there exists a constant $K$ depending only on $M$ such that
$$2 \; \leq\; \Delta f(z)\; \leq\; d + 1\; + \; Kg(z).$$ 
We leave to the reader the details of a proof using the half-space model for $\H^d$, and we refer to Lemma \ref{Laplace2} for a proof on general pinched CH manifolds.
  
\vspace{2mm}

\noindent
(b) As it will become apparent later, in dimension $d = 2$ the fact that $\Delta f(z)\sim 2$ in the neighbourhood of the geodesic line enables us to express our concentration problem in terms of the asymptotics of the first passage times for the diffusion $Y^{\nu, c, 1}$. As we said before, the spectral theory of this diffusion is somewhat simpler, because $Y^{\nu, c, 1}$ can be viewed as the square of the solution to
$$X_t \; = \; B_t \; - \; \nu \int_0^t \tanh \lpa c X_s\rpa\, ds,$$
an SDE with no more singularity at zero. In \cite{Simon:02}, the reduction to the above simple equation was already established for $d=2$, with another argument relying on Bougerol's generalized identity.
\end{remarks}

The second (Gaussian) estimate was actually already proved by Eberle - see Proposition 3.1. in \cite{Eberle:02} - for the same final purposes, though he used then the estimate in a slightly different manner - see (3.21) in \cite{Eberle:02}.

\begin{lemma} \label{E(s)} 
Let $\gamma(x,y)$ be the geodesic from $x$ to $y$ in time 1. There exist two constants $K, c > 0$ such that 
$$\PP\lcr \sup_{t\in[0,1]}\rho\lpa X^{x,y}_t, \gamma(x,y)(t)\rpa \ge u\rcr \; \le \; K e^{-c u^2}$$
for every $x, y \in \H^d_c$ and $u\ge 0$.
\end{lemma}

The third estimate is the most important one, and may have an independent interest. We present here a separated simple analytical proof for $\H^d_c$, although in the next section an even simpler probabilistic proof will be given, holding on all rank-one symmetric spaces. 

\begin{lemma} \label{heat} Let $p^{d,c}_t(y, z)$ be the heat kernel on $\H^d_c$ and $\dot\varphi(z,y)(0)$ be the unit oriented tangent vector in $z$ at the geodesic joining $z$ to $y$. Then, for every $\e \in ]0,1]$, 
$$\rho(z,y)^{-1}\grad\log p^{d,c}_t(\newdot,y)(z)\; \to \; t^{-1}\dot\varphi(z,y)(0)$$
as $\rho(z,y)\to +\infty$, uniformly on $t\in[\e,1]$ and $z,y\in \H^d_c$.
\end{lemma}

\noindent
\begin{proof} Since $p^{d,c}_t(y, z)$ only depends on $t$ and $\rho(y,z)$, we see that $\grad\log p^d_t(\newdot,y)(z)$ is parallel to $\dot\varphi(z,y)(0)$. Suppose first that $c = 1$ and set $\rho = \rho(y,z),$ $p^d_t(\rho) = p^{d,c}_t(y, z)$ and $u = \dot\varphi(z,y)(0)$ for simplicity. According to the so-called Millson's descent formula - see e.g. the fourth formula p. 5 in \cite{Anker-Ostellari:03}, we have 
$$p^{d+2}_t(\rho)\; =\; -\lpa\frac{e^{-dt}}{2\pi \sinh \rho}\rpa\frac{\partial p^d_t(\rho)}{\partial \rho}.$$
Using the closed forms of $p^d_t(\rho)$ given e.g. by (2.2) and (2.3) in \cite{Anker-Ostellari:03}, we get
\begin{eqnarray*}
\left\langle\grad\log p^d_t(\newdot,y)(z), u \right\rangle & = & \lpa 2\pi e^{dt}   \sinh \rho \rpa \frac{ p^{d+2}_t(\rho)}{p^{d}_t(\rho)}\\
& = & \lpa \sinh \rho\rpa   \frac{h^{d+2}_t(\rho)}{h^{d}_t(\rho)},
\end{eqnarray*}
with the notations
$$\lacc \begin{array}{ll} h^d_t(\rho)\; =\; \di\lpa - \f{1}{\sinh \rho} \f{\partial}{\partial \rho}\rpa^{\f{d-1}{2}} e^{-\f{\rho^2}{2t}} & \mbox{for $d$ odd,}\\
h^d_t(\rho)\; =\;  \di\int_{\rho}^{+\infty}\f{\sinh s\; ds}{\sqrt{\cosh s - \cosh \rho}}\lpa - \f{1}{\sinh s} \f{\partial}{\partial s}\rpa^{\f{d}{2}} e^{-\f{s^2}{2t}}& \mbox{for $d$ even.}
 \end{array}\right.$$
Since the involved functions are continuous with respect to $t$, by Heine's theorem it suffices to show that 
$$\lim_{\rho\to +\infty}(t\sinh \rho)h_t^{d+2}(\rho)/\rho h_t^d(\rho)\; =\; 1$$ 
for every $t > 0$. To prove this, we notice by a straightforward recurrence argument
that for every $n\in\N$ and fixed $t >0$,
\begin{equation} \label{exp}
\lpa - \f{1}{\sinh \rho} \f{\partial}{\partial \rho}\rpa^n e^{-\f{\rho^2}{2t}}\; =\; \lpa \f{\rho}{t\sinh \rho} \rpa^n \lpa 1 + O\lpa \f{1}{\rho}\rpa\rpa e^{-\f{\rho^2}{2t}},
\end{equation}
which clearly finishes the proof of the lemma when $d$ is odd. When $d$ is even, we first see that (\ref{exp}) reduces the problem to the proof of
$$\lim_{\rho\to +\infty}(t\sinh \rho){\bar h}_t^{d+2}(\rho)/\rho {\bar h}_t^d(\rho)\; =\; 1$$ 
for every $t > 0$, with the notation
$${\bar h}^d_t(\rho)\; =\;  \int_{\rho}^{+\infty}\f{\sinh s\; ds}{\sqrt{\cosh s - \cosh \rho}}\lpa \f{s}{t\sinh s} \rpa^{\f{d}{2}} e^{-\f{s^2}{2t}}.$$
This latter estimate comes now easily from the fact (whose detailed proof is left to the reader) that 
$$\lim_{\rho\to +\infty} \lpa {\bar h}^d_t(\rho) - {\bar h}^d_t(\rho + 1)\rpa/{\bar h}^d_t(\rho + 1)\; =\; +\infty.$$
This completes the proof in the case $c=1$, the case $c \neq 1$ following readily from the fact that $p^{d,c}_t(\rho) = p^{d,1}_{c^2t}(c\rho)$. 
\end{proof}

\subsection{End of the proof} We begin with the concentration around the line $\varphi(x,y)$. With the above notations, we need to prove that for every $a > 0$, 
\begin{equation}
\label{conchyp}
\lim_{s\to \infty} s^{-1}\log\PP\lcr\sup_{0\le t\le  1}Z_t(s)> a\rcr\; = \; -\KK^a_c. 
\end{equation}
Actually, from now on our method does not depend on the specific geometry of $\H^d_c$ anymore, and further on it will be readily adapted to more general manifolds, save for a comparison argument which will be detailed in the next section.

First, notice that we can replace 1 by 1/2 in the above event: once we have proved the result for $\sup_{0\le t\le 1/2}Z_t(s)$ and for every $x$ and $v$, then we can use the fact that $\lacc X_{1-t}(s), \; 0\le t\le 1\racc$ is a Brownian motion started at $y(s)$, conditioned to hit $x$ at time~$1$. In the following we will denote $X^{x,y(s)}_t$ by $X_t(s)$ for simplicity. Introducing
 $$ E(s):=\lacc\rho(X_t(s),\varphi(x,y(\infty))(st))\le s^{3/4}\ \ \hbox{for all}\ \
t\in [0,1/2]\racc,$$
we see from Lemma \ref{E(s)} that we can work on $E(s)$, i.e. it suffices to prove that
\begin{equation}
\label{conchyp2}
\lim_{s\to \infty} s^{-1}\log\PP\lcr\sup_{0\le t\le  1/2}Z_t(s)> a ; \; E(s)\rcr\; = \; -\KK^a_c. 
\end{equation}
Elementary negatively curved geometry yields the following estimates as $s\to +\infty$, uniformly on $\omega \in E(s)$ and $t\in[0,1/2]$:
$$\dot\varphi(X_t(s),y(s))(0)\;\to\;\dot\varphi(X_t(s),y(\infty))(0)\;\;\mbox{and}\;\;s^{-1}\rho(X_t(s), y(s))\;\to\;(1-t).$$
Hence, it follows from Lemma \ref{heat} that 
\begin{equation}
\label{E:9}
\lim_{s\to +\infty}s^{-1}V_t(s,X_t(s))\; = \;\dot\varphi(X_t(s),y(\infty))(0)
\end{equation}
uniformly on $\omega \in E(s)$ and $t\in[0,1/2]$. From It\^o's formula, we can now derive the following SDE for the process $Z_t(s)$: 
\begin{eqnarray*}
Z_t(s) & = & \int_0^t\left\langle df(X_u(s)),d^\n_\Ito
X_u(s)\right\rangle \; + \; 1/2\int_0^t\D f(X_u(s))\,du\\
& = & \int_0^t \left\langle \grad f(X_u(s)),A(X_u(s))\,dB_u\right\rangle \; +\; 1/2\int_0^t\D f(X_u(s))\,du\\
& + & \int_0^t \left\langle\grad f(X_u(s)), V_u(s, X_u(s)) \right\rangle\, du.
\end{eqnarray*}
First, using the formula $\di\grad f=2\sqrt{f}\f{\grad f}{\|\grad f\|}$, we can rewrite the diffusion term:
\begin{eqnarray*}
\int_0^t \left\langle \grad f(X_u(s)),A(X_u(s))\,dB_u\right\rangle  & = & \int_0^t  2 \sqrt{Z_u(s)}\,d\beta_u(s),
\end{eqnarray*}
where
$$ \beta_t(s)=\int_0^t\left\langle\f{\grad f}{\|\grad f\|}(X_u(s)),A(X_u(s))\,dB_u\right\rangle$$
is a real-valued Brownian motion for every $s > 0$. Second, we see from elementary hyperbolic geometry that
$$\left\langle\f{\grad f}{\|\grad f\|} (X_u(s)), \dot\varphi(X_u(s),y(\infty))(0)\right\rangle\; =\; -\tanh \lpa c \sqrt{Z_u(s)} \rpa.$$
Hence, it follows from (\ref{E:9}) that for every $\alpha > 0$, there exists $s_0 > 0$ such that for every $s > s_0$, 
\begin{equation}
\label{Topono1} 
 \lva s^{-1} \left\langle\f{\grad f}{\|\grad f\|} (X_u(s)), V_u(s, X_u(s)) \right\rangle + \tanh\lpa c \sqrt{Z_u(s)} \rpa \rva \; < \; \alpha,
\end{equation}
uniformly on $\omega \in E(s)$ and $u\in[0,1/2]$. Fixing now $\alpha > 0$ and  taking $s$ big enough, we deduce by comparison from (\ref{Topono1}) and Lemma \ref{Laplace} that for every $t\in [0, 1/2]$ the following a.s. inequalities hold: $Z^1_t(s) \le Z_t(s) \le Z^2_t(s)$, where $Z^1(s)$ and $Z^2(s)$ solve respectively
$$Z^1_t (s)\; =\; 2\int_0^t \sqrt{Z^1_u(s)} d\beta_u\; -\; 2s\int_0^t \sqrt{Z^1_u(s)} \lpa \tanh \lpa c \sqrt{Z^1_u(s)}\rpa + \alpha \rpa ds\; +\; t$$
and
$$Z^2_t (s)\; =\; 2\int_0^t \sqrt{Z^2_u(s)} d\beta_u\; -\; 2s\int_0^t \sqrt{Z^2_u(s)} \lpa\tanh \lpa c \sqrt{Z^2_u(s)} \rpa - \alpha \rpa ds\; +\; \f{dt}{2}.$$
This completes the proof of the concentration around $\phi(x,y)$, in letting $\alpha$ tend to 0, and using Proposition \ref{1dim2}. For the concentration around the segment $S(x,y)$, we first notice that since $S(x,y)\subset \phi(x,y)$ it is sufficient to prove 
\begin{equation}
\label{conchypbis}
\limsup_{s\to \infty} s^{-1}\log\PP\lcr\sup_{0\le t\le  1}\tilde{Z}_t(s)> a\rcr\; \le \; -\KK^a_c,
\end{equation}
where $\tilde{Z} = \tilde{f} (X)$ with $S(x,y)$ instead of $\phi(x,y)$ in the definition of $\tilde{f}$. Hence, we need to bound $\tilde{Z}_u$ a.s. from above. However, a simple picture shows that a.s. on $\{Z_u(s)\neq\tilde{Z}_u(s)\}$
$$\left\langle \f{\grad \tilde{f}}{\|\grad \tilde{f}\|}(X_u(s)), \dot\varphi(X_u(s), y(\infty))(0) \right\rangle\;\le\; -\tanh\lpa c \sqrt{\tilde{Z}_u(s)}\rpa,$$
so that by (\ref{E:9}) and Remark 2.5 (a), we see that fixing any $\alpha > 0$ and introducing
$$\tilde{Z}^2_t (s)\; =\; 2\int_0^t \sqrt{\tilde{Z}^2_u(s)} d\beta_u\; -\; 2s\int_0^t \!\!\!\sqrt{\tilde{Z}^2_u(s)} \lpa\tanh \lpa c \sqrt{\tilde{Z}^2_u(s)} \rpa - \alpha \rpa ds\; +\; \f{(d+1)t}{2},$$
for $s$ big enough and every $t\in [0, 1/2]$, the inequality $\tilde{Z}_t(s) \le \tilde{Z}^2_t(s)$ holds a.s. This allows now to obtain (\ref{conchypbis}) exactly in the same way as above.

\fin

\section{The case of rank-one noncompact symmetric spaces}\label{Rankone}\setcounter{equation}0

In this section we prove Theorem \ref{T1} on rank-one noncompact symmetric spaces, which can be viewed as a generalization of real hyperbolic spaces with pinched non constant sectional curvature. From the technical point of view, we will have to extend Lemmas \ref{Laplace}, \ref{E(s)} and \ref{heat} to this more general framework. To finish the proof, the key-argument will then consist simply in estimating the left-hand side of (\ref{Topono1}) via Alexandrov-Toponogov's comparison theorem.

\subsection{Some features of rank-one noncompact symmetric spaces}

For a complete account on the classification of such manifolds and the heat kernel thereon, we refer to Chapter X in \cite{Helgason:78} and the Anhang 4.1 in \cite{Lohoue-Rychener:82}. Let us just recall that noncompact symmetric spaces of rank-one can be divided into four families of homogeneous spaces:

\vspace{4mm}

\begin{center}
\begin{tabular}{|l|c|c|c|}
\hline
$M$ & $\alpha$ & $\beta$ & Dim $M$ \\
\hline
\hline
$\H^n(\R)\sim{\rm SO}_o(1,n)/{\rm O}(n)$ & (n-1)& 0 & n \\
\hline
$\H^n(\C)\sim{\rm SU}(1,n)/{\rm U}(n)$ & 2(n-1)& 1 & 2n \\
\hline
$\H^n(\H)\sim{\rm Sp}(1,n)/{\rm Sp}(1){\rm Sp}(n)$ & 4(n-1) & 3 & 4n \\
\hline
$\H^2(\OO)\sim F_{4}^*/{\rm Spin}(9)$ & 8 & 7 & 16 \\
\hline
\end{tabular}
\end{center}

\vspace{4mm}

\noindent
where on the three first lines we used the usual notations for classical matrix groups - see e.g. Chapter X. 2 in \cite{Helgason:78} - and ${\rm SO}_o(1,n)$ means the connected component of the identity in ${\rm SO}(1,n)$. On the last line, $F_4^*/{\rm Spin}(9)$ stands for the noncompact dual (in the sense of Chapter V. 2 in \cite{Helgason:78}) of Cayley's projective plane for octonions $\PP^2(\OO) \sim F_4/{\rm Spin}(9)$, where $F_4$ is the 52-dimensional automorphism group of ${\rm Herm}(3, \OO)$ and ${\rm Spin}(9)$ is the 36-dimensional two-fold universal covering space of ${\rm SO}(9)$ - see e.g. Section 4.1.3.5 in \cite{Berger:03} for a more complete presentation of this exceptional space. In the above table, $\alpha$ and $\beta$ are the respective multiplicities of the two generators of the root system associated with $M$. They characterize the radial part of the Laplace-Beltrami operator $\Lalb$ in $M$, whose expression is the following (see e.g. Proposition 5.26 p. 31
 3 and Formula (56) p. 315 in \cite{Helgason:84}, though we use here the probabilistic convention of Formulae (II,1) and (II,2) in \cite{Lorang-Roynette:96}):
$$\Lalb f(\rho) \; = \; \f{1}{2} f''(\rho) + \f{k}{2}(\alpha \coth k\rho + 2\beta \coth 2k\rho) f'(\rho)$$
for a positive parameter $k$. Recall that on rank-one symmetric spaces, the heat kernel $p_t(y, z)$ is a function of the sole variable $\rho(y,z)$ and that, setting $\rho = \rho(y,z)$ for simplicity, it is the fundamental solution $p_t(\rho)$ to
$$\lpa \Lalb -\partial_t\rpa f(t, \rho) \; = \; 0$$
normalized to define a probability measure on $M$. Finally, the fact that rank-one noncompact symmetric spaces have pinched negative sectional curvature $\kappa$, i.e.
\begin{equation} \label{c_2}
-c_2\; \le\; \kappa\;\le\; -c_1
\end{equation}
for some constants $c_2\ge c_1>0$, follows from the transitive action of the isometry group $G$ on $M = G/K$: on a fixed point $x\in M$, the sectional curvatures are bounded from below by smoothness of $M$, and the rank-one property yields a negative upper bound, these two bounds holding then on the whole $M$ by transitivity, since $G$ preserves the Ricci curvature tensor. The pinched property can also be seen from metric space arguments which make $M$ into a CAT(-1) space, see Theorem II.10.10 and Proposition II.10.12 in \cite{Bridson-Haefliger:99}. 

The following uniform estimate on $p_t(y,z)$, which is a direct consequence of (\ref{c_2}), Davies-Mandouvalos' estimates on real hyperbolic spaces, and a heat kernel comparison theorem - see respectively Formula (3.3) in \cite{Anker-Ostellari:03} and Theorem 4.5.2 in \cite{Hsu:02}, will be a crucial tool in extending Lemma \ref{E(s)} to rank-one symmetric spaces and more general pinched CH manifolds: there exists constants $K > 1$ and $k_2 \ge k_1 > 0,$ such that setting $\nu=(d-1)/2$ and
$$p^i_t(\rho)\; = \;t^{-d/2}(1+\rho)^\nu e^{-(k_i\rho+ \rho^2/2t)}$$
for $i = 1,2,$ there is a uniform comparison
\begin{equation}
\label{Giulini}
K^{-1} p^2_t(\rho(y,z))\; \le\; p_t(y,z)\;\le\; K p^1_t(\rho(y,z))
\end{equation}
for every $(t,y,z)\in (0,1]\times M\times M$. We stress that this uniform estimate holds on any CH manifold whose sectional curvatures satisfy (\ref{c_2}). Actually, a more precise estimate
due originally to Giulini and Mauceri - see Formula (3.1) in \cite{Anker-Ostellari:03} and the comments thereafter, holds on rank-one noncompact symmetric spaces, allowing to take $k_1 = k_2$ in (\ref{Giulini}). But we shall not use this in the sequel.

\subsection{Proof of the theorem} We first need to extend the three estimates of Section \ref{Hyperbolic}. The following extension of Lemma \ref{Laplace} holds actually on all pinched CH manifolds.
 
\begin{lemma} \label{Laplace2} Let $\phi$ be a geodesic line or segment in $M$. With the notations of Lemma \ref{Laplace}, there exists a positive constant $K$ depending only on $M$ such that
$$2 \; \leq\; \Delta f(z)\; \leq\; (d + 1) \; + \; Kg(z)$$ 
uniformly in $M$.
\end{lemma}

\begin{proof} 
We suppose first that $\phi$ is a line. Recall that
$$\Delta f(z)\; = \; 2 \;+ \; g(z) \Delta g(z)$$
so that, as in Lemma \ref{Laplace}, it remains to estimate $\Delta g(z)$. The fact that $\Delta g(z) \ge 0$ for every $z\in M$ follows from the well-known convexity property of $g$, because $g$ measures the distance to some complete convex subset of $M$ - see e.g. Corollary II.2.5 in \cite{Bridson-Haefliger:99}. To prove the other inequality, fix $z\in M$, let $z_0\in \phi$ minimize the distance between $z$ and $\phi$, and set $\rho_0$ for the distance function from $z_0$. We have
$$\Delta g(z)\; =\; \sum_{i=1}^d(g\circ\gamma_i)''(0)\;\;\;\mbox{and}\;\;\;\Delta\rho_0(z)\; =\; \sum_{i=1}^d(\rho_0\circ\gamma_i)''(0),$$
where the $\gamma_i's$ are geodesic lines such that $\gamma_i(0) = z$ and the $\dot\gamma_i(0)$'s form an orthonormal basis of $T_z M$. Besides, one can choose $\dot\gamma_1(0)$ parallel to the geodesic between $z_0$ and $z$, so that $(\rho_0\circ\gamma_1)''(0) = (g\circ\gamma_1)''(0) =0$. For every $i\ge 2$, we notice that $(g\circ\gamma_i)(0) = (\rho_0\circ\gamma_i)(0) = g(z)$, $(g\circ\gamma_i)'(0) = (\rho_0\circ\gamma_i)'(0) = 0,$ and $(g\circ\gamma_i)(t) \le (\rho_0\circ\gamma_i)(t)$ for every $t\in\R$, which entails $(g\circ\gamma_i)''(0) \le (\rho_0\circ\gamma_i)''(0)$ and finally
$$\Delta g(z)\; \leq \;\Delta \rho_0(z).$$
Now it follows from the Laplacian comparison theorem - see e.g. Theorem 3.4.2. in \cite{Hsu:02} - that 
\begin{eqnarray*}
\Delta \rho_0(z) & \le & (d-1)\sqrt{c_2}\coth\lpa \sqrt{c_2}\rho_0(z)\rpa\\
& \le & (d-1)\lpa\sqrt{c_2} +1/\rho_0(z)\rpa\; =\; (d-1)\lpa \sqrt{c_2} + 1/g(z)\rpa,
\end{eqnarray*}
where we recall that $-c_2$ is the global lower bound on the sectional curvature of $M$. This completes the proof when $\phi$ is a line. 

When $\phi$ is a segment, we can use the same arguments, save for the fact that $\Delta g$ is not continuous on the two hypersurfaces of $M$ where the distance to the geodesic line equals the distance to one of the extremities of the segment. We let the reader check by himself that this is not a major hindrance.

\end{proof}

To extend Lemma \ref{E(s)}, we will need three preparatory results. Let $x,y\in M$, $s=\rho(x,y)$, and $\varphi(x,y):\RR\to M$ be the unique geodesic parametrized by arc length such that $\varphi(x,y)(0)=x$ and
$\varphi(x,y)(s)=y$. We identify the set of  unit vectors  in $T_x M$
orthogonal to $\dot \varphi(x,y)(0)$ with the unit sphere
$S^{d-2}$. If $v\in T_x M$, we let $u\mapsto\partr0uv$ be the
parallel transport of $v$ along the geodesic $u\mapsto \varphi(x,y)(u)$. 
Recall that the Fermi coordinates of a point $z\in M$ is the unique triplet
$(u,h,\theta)\in \RR\times \RR_+\times S^{d-2}$ such that
$$z\; =\; \Phi(u,h,\theta)\; =\;\exp_{\varphi(x,y)(u)}\lpa\partr0u(h\theta)\rpa.$$

Our first preparatory result compares in Fermi coordinates the volume element of $M$ with that of $\H^d_{c_2}$. It seems to belong to the comparison folklore, even though we could not find any reference in the literature. Notice that this comparison theorem only requires that the sectional curvature of $M$ is bounded from below. 

\begin{lemma} 
\label{Fermi}
In the above Fermi coordinates, the volume element ${\rm vol}(dz)$ of $M$ is bounded from above by
$$c_2^{-(d-2)/2}\sinh^{d-2}\lpa \sqrt{c_2}h\rpa\cosh\lpa \sqrt{c_2}h\rpa dh\,du\,dV_{S^{d-2}}(\theta).$$
\end{lemma}

\begin{proof} 
Fix a point $z_0\in M$ with Fermi coordinates
$(u_0,h_0,\theta_0)$ and denote by $(e_1,\ldots,e_{d-2})$ an orthonormal
basis of the tangent space of $S^{d-2}$ at $\theta_0$. For every $h\in[0,h_0]$, consider the vector fields  $\di U(h)=\Phi_{\ast}\lpa\partial/\partial u\rpa(u_0,h,\theta_0), H(h)=\Phi_{\ast}\lpa\partial/\partial h \rpa(u_0,h,\theta_0),$ and $E_i(h)=\Phi_{\ast}\lpa e_i\rpa(u_0,h,\theta_0)$ for $i=1,\ldots, d-2$.
They are Jacobi fields along the geodesic $h\mapsto
\Phi(u_0,h,\theta_0)$, and by the Rauch and Berger comparison theorem - see e.g. Theorems 1.28 and 1.29 in \cite{Cheeger-Ebin:75} - we obtain
$$\|U(h_0)\|\;\le\;\cosh\lpa\sqrt{c_2}h_0\rpa\;\;\;\mbox{and}\;\;\;\|E_i(h_0)\|\;\le\; c_2^{-1/2}\sinh\lpa \sqrt{c_2}h_0\rpa,$$
where the respective right-hand sides correspond to $\H^d_{c_2}$. Besides, if we parametrize $S^{d-2}$ around $\theta_0$ such that $(\partial/\partial\theta^i)(\theta_0)=e_i$, then we see that he volume element at $x_0$ is smaller than or equal to
$$\|U(h_0)\|\lpa\prod_{i=1}^{d-2}\|E_i(h_0)\|\rpa \|H(h_0)\| dh\,du\,dV_{S^{d-2}}(\theta_0).$$ 
Last, since trivially $\|H(h_0)\|=1$, we obtain the desired upper bound.

\end{proof}

For the second preparatory result we will follow the proof of Lemma~3.4 in \cite{Eberle:02}, where analogous estimates are established on real hyperbolic spaces. 

\begin{lemma}\label{Eberle}
Let $X=X^{t, x,y}$ be the Brownian bridge from $x$ to $y$ in time $t$ and $\g(x,y)$ be the geodesic from $x$ to $y$ in time $1$. 
There exist three constants $a_0, K, \lambda>0$ depending only on $M$ such that  
\begin{equation} \label{a2}
\PP\lcr\rho\lpa X^{t, x,y}_{t/2}, \g(x,y)(1/2)\rpa\,\ge\, a\rcr\;\le\; K e^{-\lambda a^2/t}
\end{equation}
for every $(a,t,x,y)\in (a_0,+\infty]\times (0,1]\times M\times M$, and 
\begin{equation}\label{a4}
\PP\lcr\rho\lpa X^{t, x,y}_{t/2}, \g(x,y)(1/2)\rpa\,\ge\, a\rcr\;\le\; Kt^{-\nu/2}e^{-\lambda (a^2\wedge a^4)/t}
\end{equation}
for every $(a,t,x,y)\in (0,+\infty]\times (0,1]\times M\times M$ with $\rho(x,y) \le 1$.
\end{lemma}

\begin{proof} 
We first suppose $c_1=1$ for simplicity, the general case $c_1 > 0$ being handled in scaling the metric and the Brownian bridge. In the following calculations, the positive constants $c, K$ will depend only on $M$, but may vary from one line to another. Recalling the notation $s = \rho(x,y)$, the random variable $X^{t, x,y}_{t/2}$ has density function
$$z\;\mapsto\; q_t^{x,y}(z)\; =\;\f{p_{t/2}(x,z)p_{t/2}(y,z)}{p_t(x,y)},$$ 
an expression which can be bounded by 
$$Kt^{-d/2}(1+s)^{-\nu}(1+\rho(x,z))^{\nu}(1+\rho(y,z))^\nu e^{k_2 s-(2\rho(x,z)^2 + 2\rho(y,z)^2 - s^2)/2t}$$
in view of the estimates (\ref{Giulini}). Let now $(u,h,\theta)$ be the Fermi coordinates of $z$ with respect to $x$ and the geodesic $\phi(x,y)$. By the triangle inequality, we have 
\begin{equation}
\label{eq3}
\rho(x,z)\le |u|+h\quad\hbox{and}\quad \rho(z,y)\le |s-u|+h.
\end{equation}
On the other hand, since the triangle $x\varphi(x,y)(u)z$ has a right angle
at $\varphi(x,y)(u)$, we see that
$$\cosh\rho(x,z)\;\ge\; \cosh |u|\cosh (h)$$
by comparison with the hyperbolic case, whence
$$e^{\rho(x,z)}\;\ge\; e^{|u|}\cosh (h)\f{(1+e^{-2|u|})}{1+e^{-2\rho(x,z)}}\;\ge\;e^{|u|}\cosh(h)$$
because $|u|\le \rho(x,z)$. Writing $l(h) = \log\cosh (h)$ for concision, we finally obtain $\rho(x,z)\ge |u|+l(h)$ and, similarly, $\rho(z,y)\ge |s-u|+l(h)$. Hence
\begin{eqnarray*}
\rho(x,z)^2+\rho(z,y)^2 & \ge & \lpa|u|+l(h)\rpa^2+\lpa|s-u|+l(h)\rpa^2\\
&= &u^2+(s-u)^2+2\lpa |u|+|s-u|\rpa l(h)+2l^2(h)\\
&= &2\lpa u- s/2\rpa^2+ s^2/2+2\lpa |u|+|s-u|\rpa l(h)+2l^2(h)\\
&\ge &2\lpa u- s/2\rpa^2+ s^2/2+2sl(h)+2l^2(h).
\end{eqnarray*}
On the other hand, it follows from inequality~\eqref{eq3} that
$$\lpa1+\rho(x,z)\rpa^\nu\lpa 1+\rho(z,y)\rpa^\nu\;\le\; \lcr(1+|u|+h)(1+|s-u|+h)\rcr^\nu.$$
Plugging everything together yields
\begin{align*}
q_t^{x,y}(z)\le& \;\;
Kt^{-d/2}(1+s)^{-\nu}\lcr(1+|u|+h)(1+|s-u|+h)\rcr^\nu\\
&\qquad \times\; e^{k_2 s - 2 ((u-s/2)^2+l^2(h) + s l(h))/t}.
\end{align*}
Last, noticing that  $(1+|u|+h)(1+|s-u|+h)\le \lpa 1+ s/2+h\rpa^2$ when $u\in[0,s]$ and $\left|u-s/2\right|\ge s/2$ when $u\not\in[0,s]$, we find that
\begin{eqnarray*}
q_t^{x,y}(z) & \le & Kt^{-d/2}(1+s)^{-\nu}\lpa 1+ s/2+h\rpa^{2\nu} e^{k_2 s - ((u-s/2)^2+2l^2(h) + 2s l(h))/t }
\end{eqnarray*}
for every $t\in(0,1]$ and $x,y,z\in M$. Setting now $\mu_t^{x,y}$ for the law of $X_{t/2}^{t,x,y}$ in Fermi coordinates, it follows from Fubini's theorem and Lemma \ref{Fermi} that for every $a, b>0$
\begin{align*}
&\mu_t^{x,y}\lcr h\ge a,\ u\ge b+s/2\rcr\;\;\le\;\; Kt^{-d/2}(1+s)^{-\nu}
\int_{b+s/2}^\infty e^{-\lpa u-s/2\rpa^2/t}\, du\\ 
&\times\int_a^\infty \lpa 1+s/2+h\rpa^{2\nu}\sinh^{d-2}\lpa\sqrt{c_2}h\rpa\cosh\lpa \sqrt{c_2}h \rpa e^{k_2 s -2(s l(h) + l^2(h))/t}\, dh.
\end{align*}
The first integral on the right hand side can be estimated followingly:
\begin{equation}
\label{eq5}
\int_{b+s/2} e^{-\lpa u-s/2\rpa^2/t}\, du\;=\;\int_b^\infty
e^{- v^2/t}\, dv\;\le\; Kt^{1/2}e^{-b^2/2t}.
\end{equation}
For the second integral, notice first that there exists $\l_1>0$
such that $l(h)\ge \l_1 h$ for every $h\ge 1$. Consequently, recalling that $t\in (0,1]$, the second integral can be bounded by 
$$K\int_a^\infty e^{- (3\l_1^2h^2)/2t}\, dh\;\le\; K t^{1/2}e^{-\l_1^2 a^2/t}$$
when $a\ge a_0 = k_2/2\l_1 $. Together with~\eqref{eq5}, this latter bound yields 
\begin{equation} \label{eq6}
\mu_t^{x,y}\lcr h\ge a,\ u - s/2 \ge b\rcr\;\le\; Kt^{-\nu}e^{-(b^2/2 + \l_1^2a^2)/t}
\end{equation}
and by symmetry, this estimate is also valid for $\di \mu_t^{x,y}\lcr h\ge a,\ u -s /2 \le -b\rcr$.
We can now establish the estimate (\ref{a2}), since
\begin{eqnarray*}
\PP\lcr\rho\lpa X^{t, x,y}_{t/2}, \g(x,y)(1/2)\rpa\,\ge\, a\rcr &\le & \mu_t^{x,y}\lcr h\ge a/2\rcr+\mu_t^{x,y}\lcr\left|u-s/2\right|\ge a/2\rcr\\
&\le & Kt^{-\nu}e^{-ca^2/t} \\
&\le & K e^{- ca^2/t} 
\end{eqnarray*}
for some constants $c, K >0$ and every $(a,t,x,y)\in [a_0, +\infty)\times (0,1]\times M\times M$. 

On the other hand, we can choose $\l_2>0$ such that $l(h)\ge \l_2 h^2$ for every $h\in [0,a_0]$. Hence, when $a\in(0,a_0]$ and $s \le 1$, the second integral from $a$ to $a_0$ is bounded by 
$$K\int_a^{a_0} e^{- 2\l^2_2h^4/t}h^{d-2}\, dh\; \le\; Kt^{\nu/2}\int_{a\l_2^{1/2}t^{-1/4}}^{\infty}e^{-2r^4}r^{d-2}\,dr\;\le\; Kt^{\nu/2}e^{-(\l_2^2a^4)/t}.$$
Together with~\eqref{eq5} and \eqref{eq6}, this entails
$$\mu_t^{x,y}\lcr  h\ge a,\ u - s/2 \ge b\rcr\;\le\; Kt^{-\nu/2}e^{-t^{-1}(b^2/2 + \l_2^2a^4)/t}$$
for every $(a,t,x,y)\in (0,a_0]\times (0,1]\times M\times M$ with $\rho(x,y) \le 1$, and we can then obtain the estimate (\ref{a4}) similarly as above.

\end{proof}

\begin{remark} \label{precise} When $M$ is a rank-one noncompact symmetric space, the fact that we can take $k_1 = k_2$ in (\ref{Giulini})  yields a better result:
$$\PP\lcr\rho\lpa X^{t, x,y}_{t/2}, \g(x,y)(1/2)\rpa\,\ge\, a\rcr\;\le\; Kt^{-\nu/2}e^{-\lambda (a^2 \wedge a^4)/t}$$
for every $(a,t,x,y)\in (0,+\infty]\times (0,1]\times M\times M$ - this is actually inequality (3.6) in \cite{Eberle:02}. However, as we said before, we shall not need this in the sequel.  

\end{remark}

For every integer $n$ and for $i = 1 \ldots 2^n$, set $t^i_n = i 2^{-n}$ and consider the event
\begin{eqnarray*}
\Lambda^{x,y}_n & = & \bigcup_{p\ge n} \lacc \sup_{1 \le i\le 2^p}\rho\lpa X^{x,y}_{t^i_p}, X^{x,y}_{t^{i-1}_p}\rpa \ge 1\racc.
\end{eqnarray*}
Our last preparatory result is an estimate on the uniform continuity of $X^{x,y}$ when $\rho(x,y) \to +\infty$.
 
\begin{lemma} 
\label{n(s)}
Setting $s = \rho(x, y)$ and $n(s) = 2 + [2\log s/\log 2]$, there exist constants $c, K > 0$ such that
$$\PP\lcr \Lambda^{x,y}_{n(s)}\rcr \; \le \; K e^{-c s^2}$$
for every $x, y \in M$.
\end{lemma}

\begin{proof} For every integer $p$ and $i = 2 \ldots 2^{p}-1$ we get from the estimates (\ref{Giulini})  
\begin{eqnarray*}
\PP\lcr \rho\lpa X^{x,y}_{t^i_p}, X^{x,y}_{t^{i-1}_p}\rpa \ge 1\rcr & = & \int_{\rho(z,t) \ge 1} \frac{p_{t^{i-1}_p}(x,z)p_{2^{-p}}^{}(z,t)p_{1-t^i_p}(t,y)}{p_1(x,y)}{\rm vol}(dz){\rm vol}(dt)\\
& \le &  \!\!\!\!\sup_{\rho(z,t) \ge 1} p_{2^{-p}}^{}(z,t)\int \frac{p_{t^{i-1}_p}(x,z)p_{1-t^i_p}(t,y)}{p_1(x,y)}{\rm vol}(dz){\rm vol}(dt)\\
& \le & \!\!\!\! \sup_{\rho(z,t) \ge 1}  p_{2^{-p}}^{}(z,t) / p_1(x,y)\\
& \le & K2^{pd/2} e^{(s^2 - 2^p)/2}, 
\end{eqnarray*}
and it is easy to see that the same inequality holds when $i = 1$ or $2^p$. Hence, for every $p \ge n(s)$, we have
\begin{eqnarray*}
\PP\lcr \sup_{1 \le i\le 2^p} \rho\lpa X^{x,y}_{t^i_p}, X^{x,y}_{t^{i-1}_p}\rpa \ge 1\rcr & \le &  K2^{p(d/2 + 1)} e^{-2^p/4}\; \le\; K e^{-2^p/5},
\end{eqnarray*}
which readily entails
$$\PP\lcr \Lambda^{x,y}_{n(s)}\rcr\; \le\; K e^{-2^{n(s)}/5} \; \le \; K e^{-c s^2}$$
for some constants $c, K > 0$ independent of $x, y\in M$.

\end{proof}

The following proposition yields a weak extension of Lemma \ref{E(s)} on pinched CH manifolds. Actually on rank-one symmetric spaces the exact statement of Lemma \ref{E(s)} could be transfered verbatim, because of Remark \ref{precise}. However, on pinched CH manifolds the extension takes the form of a limsup theorem because the (optimal) inequalities (\ref{Giulini}) are not precise enough to allow a uniform estimate. Nevertheless, as we see from the end of the proof in Section 2, the result will be sufficient for our purposes. The proof mimics that of Proposition~3.1 in \cite{Eberle:02}, save for the use of Lemma \ref{n(s)}.

\begin{prop} \label{E2(s)} Fixing $x\in M$, $v\in T_x M$ unitary and setting $y = y(s) = \exp_x(sv)$, there exists two constants $c, K > 0$ independent of $s, v$ such that 
$$\PP\lcr \sup_{t\in[0,1]}\rho\lpa X^{x,y}_t, \gamma(x,y)(t)\rpa \ge s^{3/4}\rcr \; \le \; K e^{-c s^{5/4}}.$$
\end{prop}

\begin{proof} Clearly, we can suppose that $s \ge 1$ and again, in the following calculations the positive constants $c, K$ will depend only on $M$ but may vary from one line to another. Let $\SP=\SP_{x,y}$ be  the set of continuous paths $\om  : [0,1]\to M$ satisfying $\om(0)=x$ and $\om(1)=y$, and $\PP^{x,y}$ be the law of the Brownian bridge on $\SP$. For $k\in \NN$ and $\om\in\SP$, set
$$M_k(\om)=\max_{0\le i\le 2^k}\rho\lpa\om(t^i_k),\g(x,y)(t^i_k)\rpa$$
and 
$$N_k(\om)=\max_{0\le i< 2^k}\rho\lpa\om\lpa t^{2i+1}_{k+1}\rpa,\g(\om(t^i_k),\om(t^{i+1}_k)(1/2)\rpa.$$
By convexity of the distance fonction $\rho$ on $M\times M$ - see e.g. Proposition II.2.2 in \cite{Bridson-Haefliger:99}, we see that $\di M_{k+1}(\om)\le M_k(\om)+N_k(\om)$. Since $M_0(\om)=0$ and by continuity of the path $\om$, we have
\begin{equation}
\label{eq8}
\sup_{t\in[0,1]}\rho\lpa \om(t),
\g(x,y)(t)\rpa\;\le\;\sum_{j=0}^{+\infty} N_j(\om).
\end{equation}
Setting 
$$\Omega^{x,y}\; =\; \lacc\sup_{t\in[0,1]}\rho(X^{x,y}_t,\g(x,y)(t))\ge s^{3/4}\racc,$$
we deduce from Lemma \ref{n(s)} that it suffices to prove that
\begin{equation} \label{oxy}
\PP\lcr \Omega^{x,y} \cap \lpa\Lambda^{x,y}_{n(s)}\rpa^c\rcr\; \le \; K e^{-c s^{5/4}}.
\end{equation}
Introducing the event
$$\Lambda_n \; = \; \bigcup_{p\ge n} \lacc \sup_{1 \le i\le 2^p}\rho\lpa \omega(t^i_p), \omega(t^{i-1}_p)\rpa \ge 1\racc,$$
we see from (\ref{eq8}) that 
\begin{eqnarray*}
\PP\lcr \Omega^{x,y}_{\delta}\cap \lpa\Lambda^{x,y}_{n(s)}\rpa^c\rcr &\le & \sum_{j=0}^{+\infty}\;\PP^{x,y}\lcr\lacc N_j\ge c 2^{-j/16}s^{3/4}\racc \cap\Lambda_{n(s)}^c \rcr.
\end{eqnarray*}
On the one hand, from estimate (\ref{a2}) and the fact that $2^{-n(s)} = c s^{-2}$, for $s$ big enough we have
\begin{eqnarray*}
\sum_{j=0}^{n(s)}\;\PP^{x,y}\lcr N_j\ge c 2^{-j/16}s^{3/4} \rcr
&\le & \sum_{j=0}^{n(s)}\;\PP^{x,y}\lcr N_j\ge c s^{5/8} \rcr\\
& \le & \sum_{j=0}^{n(s)} \;2^j\!\! \sup_{z,t \in M}\mu_{2^{-j}}^{z,t}\lcr \rho\lpa\cdot,\g(z,t)(1/2)\rpa\ge c s^{5/8}\rcr\\
&\le & K s^2 n(s) e^{- cs^{5/4}}\;\le\; K e^{-cs^{5/4}}.
\end{eqnarray*}
On the other hand, from estimate (\ref{a4}),
\begin{align*}
&\sum_{j= n(s)}^{+\infty}\PP^{x,y}\lcr\lacc N_j\ge c 2^{-j/16}s^{3/4}\racc \cap\Lambda_{n(s)}^c \rcr \\
&\qquad\le\;\;\sum_{j= n(s)}^{+\infty}\; 2^j\!\!\!\! \sup_{\tiny\begin{array}{c}z,t \in M\\ \rho(z,t) \le 1\end{array}}\mu_{2^{-j}}^{z,t}\lcr \rho\lpa\cdot,\g(z,t)(1/2)\rpa\ge c 2^{-j/16} s^{3/4}\rcr\\
&\qquad\le\;\; K \sum_{j= n(s)}^{+\infty}2^{j\lpa\nu/2+1\rpa} e^{-c(2^{3j/4} s^{3})\wedge (2^{7j/8} s^{3/2})}\; \le \; K e^{-c2^{3n(s)/4} s^{3/2}}\; \le \; Ke^{-cs^2}.
\end{align*}
Gluing these two latter estimates together yields (\ref{oxy}), and completes the proof of the proposition.

\end{proof}

Finally, the next proposition extends Lemma \ref{heat} to rank one noncompact symmetric spaces. Its proof relies on a nice probabilistic representation of the heat kernel on the latter, which is due to Lorang and Roynette \cite{Lorang-Roynette:96}. We notice that their closed formula carries over more general Sturm-Liouville operators on $\R^+$, and hence allows to consider e.g. radially symmetric manifolds whose sectional curvature is constant and negative at infinity. Nevertheless, the extension of this estimate to general pinched CH manifolds - see Section \ref{CH} - will be a more difficult task, requiring sophisticated probabilistic tools.
 
\begin{prop} \label{heat2} Let $p_t(y, z)$ be the heat kernel on $M$ and $\dot\varphi(z,y)(0)$ be the unit oriented tangent vector in $z$ at the geodesic joining $z$ to $y$. Then, for every $\e \in ]0,1]$, 
$$\rho(z,y)^{-1}\grad\log p_t(\newdot,y)(z)\; \to \; t^{-1}\dot\varphi(z,y)(0)$$
as $\rho(z,y)\to +\infty$, uniformly on $t\in[\e,1]$ and $z,y\in M$.
\end{prop}

\begin{proof}  Set $u = \dot\varphi(z,y)(0)$, $\rho = \rho(y,z)$ and suppose first $d > 2, k=1$. According to formula (II,21) in \cite{Lorang-Roynette:96}, we can express
\begin{eqnarray*}
p^1_t(\rho) & = & t^{-d/2}(1+\rho)^{(d-1)/2}\exp-\lpa(\alpha + 2\beta)^2 t/8 + (\alpha + 2\beta)\rho/2 + \rho^2/2t\rpa\Phi_t(\rho)
\end{eqnarray*}
where $\lacc X^{t, \rho}_s, \; s\in [0,t]\racc$ is the $d$-dimensional Bessel bridge from $\rho$ to 0 in time t, $\tilde{l}$ is the real function given by (II, 24) in \cite{Lorang-Roynette:96}, and
$$\Phi_t(\rho)\; =\; \E\lcr \exp - \lpa\int_0^t\tilde{l}(X^{t, \rho}_s)\, ds\rpa\rcr.$$
Hence, we just need to prove that $\Phi_t(\rho)'/\rho\Phi_t(\rho)\to 0$ when $\rho\to +\infty$ for a fixed $t > 0$. Since $\tilde{l}$ and $\tilde{l}'$ are bounded functions, by dominated convergence this clearly amounts to prove that
$$\lim_{\rho \to +\infty} \rho^{-1} Y^{t,\rho}_s\; = \; 0\;\;\;\;\mbox{a.s.}$$
where $Y^{t,\rho}_s = D_{\rho} X^{t,\rho}_s$ solves the random ODE (see the Appendice in \cite{Lorang-Roynette:96})
\begin{eqnarray*}
Y^{t,\rho}_s & = & 1 - \int_0^s Y^{t,\rho}_u \lpa \f{(d-1)}{2} (X^{t,\rho}_u)^{-2}+ \f{1}{(s-u)}\rpa \, du\\
& = & \exp \lcr - \int_0^s  \lpa \f{(d-1)}{2}(X^{t,\rho}_u)^{-2}+ \f{1}{(s-u)}\rpa \, du\rcr.
\end{eqnarray*}
It is clear that $0 < Y^{t,\rho}_s < 1$ a.s. and the proof is complete in the case $d > 2, k=1$. Again, the case $d > 2, k\neq 1$ follows readily from the fact that $p^k_t(\rho) = p^1_{k^2t}(k\rho)$, with the obvious notations. Finally, the case $d = 2$ consists only, up to isomorphism, in $\H^2_{c}(\R)$, and this case was already treated in Lemma \ref{heat}. 

\end{proof}

\begin{remark} In some sense, the above lemma says that the leading term for $p_t(y,z)$ when $\rho(y,z)$ tends to $+\infty$ with $t$ bounded away from 0 and $+\infty$, looks like $e^{-\rho^2(y,z)/2t}$. Notice that the pinched negative sectional curvature of rank-one symmetric spaces does not seem to play a r\^ole for this estimate, since on $\R^d$ there is of course an equality
$$\rho(z,y)^{-1}\grad\log p^d_t(\newdot,y)(z)\; = \; t^{-1}\dot\varphi(z,y)(0)$$
for every $y, z \in\R^d$ and $t > 0$. Even though this is no more relevant to our purposes - see the counterexample below, we believe that the estimate of Lemma \ref{heat2} also holds in the higher rank case. However, this task probably demands a more detailed analysis, since here $p_t(y, z)$ is not a function of one variable anymore. 

\end{remark}

\vspace{2mm}

\noindent
{\bf End of the proof}. Because of Lemma \ref{Laplace2}, Propositions \ref{E2(s)} and \ref{heat2}, we can actually reason as in Section 2.3 almost literally. The only point to change is (\ref{Topono1}), which does not hold anymore because the sectional curvature is not constant in general. However, Alexandrov-Toponogov's theorem - see e.g. Theorem 73 in \cite{Berger:03} - allows to compare
$$\left\langle \f{\grad f}{\|\grad f\|} (X_u(r)),\dot\varphi(X_u(r),y(\infty))(0)\right\rangle$$ 
with the same quantities on manifolds with constant curvature $-c_1$ and $-c_2$ respectively. For every $\alpha > 0$, this yields the existence of $s_0 > 0$  such that for every $s > s_0$, 
\begin{eqnarray*}
s^{-1}\left\langle\f{\grad f}{\|\grad f\|} (X_u(s)), V_u(s, X_u(s)) \right\rangle & \ge & - \tanh \lpa c_2 \sqrt{Z_u(s)} \rpa + \alpha
\end{eqnarray*}
and 
\begin{eqnarray*}
s^{-1}\left\langle\f{\grad f}{\|\grad f\|} (X_u(s)), V_u(s, X_u(s)) \right\rangle & \le & - \tanh \lpa c_1 \sqrt{Z_u(s)} \rpa -\alpha
\end{eqnarray*}
uniformly on $\omega \in E(s)$ and $u\in[0,1/2]$. We can then finish the proof exactly as in Section 2.3.

\fin

\begin{remark} The proof of Lemma \ref{heat2} shows actually a uniform speed of convergence in O($\rho^{-1}$). Since then $s\alpha$ is bounded, this makes it possible to finish the proof of Theorem \ref{T1} without Proposition \ref{1dim2}, in using Proposition \ref{1dim1} and the Cameron-Martin transformation.
\end{remark}

\subsection{A counterexample in rank two} The fact that Euclidean Brownian bridges $\lacc B^{x,y}_t, \; 0\leq t \leq 1\racc$, with one extremity $x$ far away from the other $y$, do not concentrate on geodesic lines, follows easily from the a.s. representation 
$$B^{x,y}_t\; = \; (1-t)x + ty + B_t -tB_1$$
for every $t \in [0,1]$, where $\lacc B_t, \; t \geq 0\racc$ is a Brownian motion starting from 0. In this paragraph, we would like to point out that our concentration problem for the Brownian bridge is also irrelevant on higher rank noncompact symmetric spaces, in describing an elementary counterexample on the bidisk. This space is the cartesian product $\H = \H_1\times\H_2$, where $\H_1$ and $\H_2$ are two copies of the real hyperbolic plane with sectional curvature -1, and it is the simplest example of a rank-two non compact symmetric space. We endow it with the Riemannian distance 
$$d(u,v) = \sqrt{\rho^2_1(u_1, v_1) + \rho^2_2(u_2, v_2)},$$ 
the notations being obvious. Taking rectangular coordinates on the half-space model and considering the geodesic line $\Lambda = \lacc (0,1,0, e^t),\; t\in\R\racc$, we see that the distance from $\Lambda$ of the Brownian bridge $\lacc (X^1_u, X^2_u), \; u\in [0,1]\racc$ between $(0,1,0,1)$ and $(0,1,0,e^s)$ is greater than $\rho_1(X^1_u, (0,1))$, which is a random variable independent of the parameter $s$. 

Notice that in this counterexample, we took a point at infinity in $\H_1\times\partial\H_2$. But we stress that the Brownian bridge does not concentrate either around geodesics, at least exponentially, when taking a point at infinity in $\partial \H_1\times\partial\H_2$. To see this, consider the geodesic line $\Lambda = \lacc (0,e^t,0, e^t),\; t\in\R\racc$, let $\PP^s$ be the law of the Brownian bridge between $(0,1,0,1)$ and $(0,e^s,0,e^s)$, and $\lacc X_u = (X^1_u, X^2_u), \; u\in [0,1]\racc$ be the coordinate process on $\H$. Setting $i$ for the identity map between $\H_1$ and $\H_2$ and $Y^2_u = i(X^1_u)$, it follows from the definition of $d$ that for every $t\in[0,1]$ and $a >0$, 
\begin{equation}
\label{omega}
\Omega^a_t \; =\; \lacc \rho_2\lpa X^2_t, Y^2_t\rpa > 2a\racc\; \subset\; \lacc \sup_{u\in [0,1]} d(X_u, \Lambda) > a\racc.
\end{equation}
If now $\tilde{\PP}^s$ stands for the law of the Brownian motion $B$ on $\H$ starting from $(0,1,0,1)$ and conditioned to go to $(0,\infty,0,\infty)$ with speed $s$, a straightforward computation using the representation of $B$ as an exponential functional of Euclidean Brownian motion - see Part 2 in \cite{Simon:02} for details - yields
$$\tilde{\PP}^s \lcr \Omega^a_t \rcr \; \ge \; \PP\lcr \lva W^1_t - W^2_t \rva > 2a\rcr,$$
where $W^1$ and $W^2$ are two independent linear Brownian motions. Since the right-hand side does not depend on $s$, it remains to compare $\PP^s \lcr \Omega^a_t \rcr$ and $\tilde{\PP}^s \lcr \Omega^a_t \rcr$, and this can be done analogously as in \cite{Simon:02} pp. 1986-87: for every $\alpha > 0$, we can prove that 
$$\lim_{s\to +\infty} s^{\alpha}\log \PP^s\lcr \Omega^a_t\rcr \; =\; 0,$$
which proves that there is no exponential concentration, in view of (\ref{omega}). Still, one may ask if a polynomial concentration occurs when taking a point at infinity in $\partial\H_1\times\partial\H_2$. This is partly motivated by the fact that $\H_1\times\partial\H_2$ and $\partial\H_1\times\partial\H_2$ play actually an entirely different r\^ole as a subset of the topological boundary of $\H$ - we refer to \cite{Carlo-Maria:93} for much more on this topic.

\section{The case of pinched Cartan-Hadamard manifolds}\label{CH}\setcounter{equation}0

We begin with a - non-uniform - extension of Proposition \ref{heat2} to pinched CH manifolds, a result which one may find interesting in its own right:

\begin{thm} \label{heat3} Let $p_t(x, y)$ be the heat kernel on $M$. Fixing $z\in M$, for every $\e \in ]0,1]$
$$\rho(z,y)^{-1}\grad\log p_t(\newdot,y)(z)\; \to \; t^{-1}\dot\varphi(z,y)(0)$$
when $\rho(z,y)\to +\infty$, uniformly on $t\in[\e,1]$ and $y\in M$.
\end{thm}

The proof of this theorem is probabilistic and entirely independent of the preceding sections. It relies on a Bismut-type formula yielding a representation of the logarithmic derivative of the heat kernel in terms of the Brownian bridge \cite{Thalmaier:97}, and suitable large deviations estimates relying on Varadhan's lemma. Using analogous arguments, recall that Bismut, motivated by Brownian holonomy and probabilistic index theory  - see Theorem 3.8 in \cite{Bismut:84} - had proved that on a compact manifold
$$t\grad\log p_t(\newdot,y)(z)\; \to \; \rho(z,y)\dot\varphi(z,y)(0)$$
when $t\to 0$, provided that $y$ and $z$ are not in each other's cut locus. We also refer to Theorem 2.5 in \cite{Norris:93} for an extension of this limit theorem to the successive derivatives of the heat kernel.

We will need a preparatory result, and for this we fix some notations. Let $z\in M$, $v\in T_z M$ unitary, and for every $s\ge 0$, set $y = y(s) = \exp_z(sv)$, so that $s = \rho(z,y)$ and $v = \dot\varphi(z,y)(0)$. Let $\lacc Y_u,\; u\ge 0\racc$ be a Brownian motion in $M$ started at $z$. Fix $\b\in(0, 1]$, set $r = 1/(\b s)$ and $\lacc Y^r_u = Y_{ru}, \; u\ge 0\racc$ for the Brownian motion with speed $r$. Let $\lacc \FF^r_u,\; u\ge 0\racc$ be its canonical completed filtration and consider the $\FF^r_u$-stopping times
$$T^{\l}_r \; = \; \inf\{t>0,\; \rho(z, Y^r_t) \ge \l\}$$
for every $\l > 0$. The following lemma will be useful for a crucial localisation procedure during the proof of Theorem \ref{heat3}.

\begin{lemma} 
\label{cLambda} Fix $\e > 0$. For every $c > 0$, there exists $\l>0$ such that
$$\PP\lcr T^{\l}_r<1\;\vert\; Y^r_{t/r}=y \rcr \; \le \; e^{-c s}$$
uniformly on $z\in M,$ $v\in T_z M$ and $t \in (\e, 1]$.
\end{lemma}

\begin{proof} It follows from the inhomogeneous Markov property that 
\begin{eqnarray*}
\PP\lcr T^{\l}_r<1\;\vert\; Y^r_{t/r}=y \rcr & = & \E\lcr \Un_{\lacc T^{\l}_r<1\racc}\frac{p_{t-(rT^{\l}_r)}(Y^r_{T^{\l}_r},y)}{p_t(z,y)}\rcr\\
& \le &  \PP\lcr T^{\l}_r<1\rcr\!\!\!\sup_{\tiny\begin{array}{c}\rho(z,x) \le \l\\ u \le r\end{array}}\lpa\f{p_{t-u}(x,y)}{p_t(z,y)}\rpa
\end{eqnarray*}
On the one hand, it follows from (\ref{Giulini}) that for every $u\le r$, every $s \ge (2/ \e)\vee\lambda$ and every $x$ such that $\rho(z,x) \le \l$,
\begin{eqnarray}
\label{inhom}
\f{p_{t-u}(x,y)}{p_t(z,y)}&\le & K e^{(s^2-\rho^2(x,y))/2t  + k_2 s}\;\le\; Ke^{k(1+\lambda)s}
\end{eqnarray}
for some positive constants $k, K$ independent of $z\in M,$ $v\in T_z M,$ $t \in (\e, 1]$, $\lambda$ and $s$. On the other hand,
\begin{eqnarray*}
\PP\lcr T^{\l}_r<1\rcr & \le & \PP\lcr \sup_{u\le 1} \rho(z, Y^r_u) \ge \l\rcr\; =\; \PP\lcr \sup_{u\le r} \rho(z, Y_u) \ge \l\rcr
\end{eqnarray*}
and since the process $u \mapsto \rho(z, Y_u)$ lies a.s. between two Bessel processes whose parameters depend only on $c_1, c_2$ - see e.g. Corollary 3.4.4 and Theorem 3.5.1 in \cite{Hsu:02}, an immediate scaling argument yields
\begin{eqnarray}
\label{bess}
\PP\lcr \sup_{u\le r} \rho(z, Y_u) \ge \l\rcr & \le & e^{-k'\l^2\b s}
\end{eqnarray}
for some constant $k'$ depending only on $c_1, c_2$. Putting (\ref{inhom}) and (\ref{bess}) together completes the proof of the lemma.

\end{proof}

\noindent
{\bf Proof of Theorem \ref{heat3}.} We will use the same notations as above and, for concision, in the following we will write "uniformly" for "uniformly on $t\in[\e,1]$ and $y\in M$". Let $\lacc h_u, \; 0\le u\le 1\racc$ be an a.s. differentiable, $\FF^r_u$-adapted nonincreasing process with values in $[0,1]$ such that $h_0=1$, $h_1=0$, $\dot h_u=-1$ when $u \le T^{\l}_r \wedge 1$, 
$h_u= 0$ when $t \ge T^{\l + 1}_r$, and $\int_0^1\dot h_u^2\,du\in L^2$ - we refer to \cite{Thalmaier-Wang:98} for the construction of $h$. From Corollary~2.5 and Formula~(6.7) in~\cite{Thalmaier:97} together with a Brownian scaling argument, we have for every $r > 0$
\begin{eqnarray*}
r\grad\log p(t,\newdot, y)(z) & = & -\E\lcr\int_0^1\dot h_w(\Theta_{0,w}^r)^\ast d_\Ito^\n Y^r_w\;\vert\; Y^r_{t/r}=y\rcr\nonumber\\
& = &-\E\lcr\lpa\int_0^1\dot h_w(\Theta_{0,w}^r)^\ast d_\Ito^\n
  Y^r_w\rpa\f{p_{t-r}(Y^r_1,y)}{p_t(z,y)}\rcr,
\end{eqnarray*}
where $\Theta_{0,u}^r : T_zM\to T_{Y^r_u}M$ is the so-called deformed parallel translation along $Y^r$: $\Theta_{0,0}^r=\Id_{T_zM}$ and for every $g\in T_zM$, the process $\lacc \Theta_{0,u}^r(g), \; u\ge 0\racc$ satisfies the random covariant equation
$$D\Theta_{0,u}^r(g)\;=\; -\f{r}{2}\Ric^\sharp\lpa\Theta_{0,u}^r(g)\rpa du.$$ 
Using the Markov property and the fact that $h_u= 0$ when $t \ge T^{\l + 1}_r$, we get the following decomposition
\begin{align}\label{Tlr}
& s^{-1}\grad\log p(t,\newdot, y)(z)\; = \;  \b \E\lcr \Un_{\lacc T^{\l}_r\ge 1\racc}\lpa\int_0^1(\Theta_{0,w}^r)^\ast d_\Ito^\n Y^r_w\rpa\f{p_{t- r}(Y^r_1,y)}{p_t(z,y)}\rcr\nonumber\\
 & \;-\; \b\E\lcr\Un_{\lacc T^{\l}_r<1\racc}\lpa\int_0^{S^{\l}_r}\!\!\!\dot h_w(\Theta_{0,w}^r)^\ast d_\Ito^\n Y^r_w\rpa\f{p_{t -(rS^{\l}_r)}(Y^r_{S^{\l}_r},y)}{p_t(z,y)}\rcr
\end{align}
for every $\b > 0$, having set $S^{\l}_r = 1\wedge T^{\l + 1}_r$ for simplicity. By the Cauchy-Schwarz inequality and since $\b \le 1$, the second summand on the right-hand side in (\ref{Tlr}) is smaller than
$$\PP\lcr T^{\l}_r<1\;\vert\; Y_{t/r}=y\rcr^{1/2}\!\!\!\!\sup_{\tiny\begin{array}{c}\rho(z,x) \le \l +1\\ u \le r\end{array}}\!\!\!\!\lpa\f{p_{t-u}(x,y)}{p_t(z,y)}\rpa\left\|\int_0^1\dot h_w(\Theta_{0,w}^r)^\ast d_\Ito^\n Y^r_w\right\|_2.$$
Since by assumption $\int_0^1\dot h_u^2\,du\in L^2$ - with a norm independent of $t, z, v$ and $s$, and since a.s.  $\vert\vert (\Theta_{0,w}^r)^\ast\vert\vert \le e^{(d-1)c_2}$ for every $w\le 1$ - see Inequality (7.4) in \cite{Thalmaier:97},  we get the upper bound
$$K\PP\lcr T^{\l}_r<1\;\vert\; Y_{t/r}=y\rcr^{1/2}\!\!\!\!\!\!\!\!\!\sup_{\tiny\begin{array}{c}\rho(z,x) \le \l +1\\ u \le r\end{array}}\!\!\!\!\!\lpa\f{p_{t-u}(x,y)}{p_t(z,y)}\rpa$$
which, reasoning as in Lemma \ref{cLambda}, is smaller than
$$K e^{cs} \PP\lcr T^{\l}_r<1\;\vert\; Y_{t/r}=y\rcr^{1/2}$$
for some positive constants $c, K$ independent of $t, z, v$ and $s$. But according to Lemma \ref{cLambda}, this last expression tends to 0 uniformly when $s\to + \infty$, provided that $\l$ is big enough. Hence, we see that the second summand in the decomposition (\ref{Tlr}) is negligible and, fixing $\l > 0$ big enough once and for all,  it remains to prove that 
\begin{eqnarray*}
\b \E\lcr \Un_{\lacc T^{\l}_r\ge 1\racc}\lpa\int_0^1(\Theta_{0,w}^r)^\ast d_\Ito^\n Y^r_w\rpa\f{p_{t- r}(Y^r_1,y)}{p_t(z,y)}\rcr & \rightarrow & t^{-1} v
\end{eqnarray*}
uniformly when $s\to + \infty$. Since $\PP\lcr T^{\l}_r\ge 1\;\vert\; Y_{t/r}=y\rcr \to 1$ uniformly when $s\to + \infty$ - again according to Lemma \ref{cLambda}, this amounts to prove that
\begin{eqnarray}\label{finalresult}
\lim_{\b\to 0}\lpa\limsup_{s\to +\infty}
\E\lcr \Un_{\lacc T^{\l}_r\ge 1\racc} \left\| \Phi^r_1 -v\right\| \f{p_{t- r}(Y^r_1,y)}{p_t(z,y)}\rcr\rpa & = & 0
\end{eqnarray}
uniformly, where $\lacc \Phi_u^r,\; u\ge 0\racc$ is the $T_z M$ valued process defined by
$$\Phi_u^r\; =\; \b t\int_0^u(\Theta_{0,w}^r)^{\ast}\lpa d_\Ito^\n Y_w^r\rpa$$
for every $u\ge 0$.  Setting 
$$A_{\a}\;=\;\lacc\phi\in C([0,1],\ T_zM),\ \left\|\phi(1)-v\right\|\ge\a\racc$$ 
for every $\a > 0$, and writing $A_{\a}^r = \lacc T^{\l}_r\ge 1\racc\cap\lacc (\Phi^r)^{-1}(A_{\a})\racc$ for concision, we have
\begin{eqnarray}\label{expect}
\E\lcr \Un_{\lacc T^{\l}_r\ge 1\racc} \left\| \Phi^r_1 -v\right\| \f{p_{t- r}(Y^r_1,y)}{p_t(z,y)}\rcr & = &\int_0^{\infty}\E\lcr  \Un_{A_{\a}^r}\f{p_{t-r}(Y_1^r,y)}{p_t(z,y)}\rcr d\a,
\end{eqnarray}
so that we need an upper bound on
$$\E\lcr  \Un_{A_{\a}^r}\f{p_{t-r}(Y_1^r,y)}{p_t(z,y)}\rcr.$$
It follows from (\ref{Giulini}) that on $\lacc T^{\l}_r\ge 1\racc$ there is a constant $K$ independent of $t, z, v$ and $ s$, such that 
\begin{eqnarray*}
\f{p_{t-r}(Y_1^r,y)}{p_t(z,y)}&\le & K e^{(s^2-\rho^2(Y^r_1,y))/2t  -rs/2t(t-r) + k_2 s}\\
&\le & K e^{k_2s - s/2\b t^2 - s\psi^{}_{\infty}(Y_1^r)/t},
\end{eqnarray*}
where $\psi^{}_{\infty}$ is the Busemann function associated to $y(\infty)$ and vanishing at $z$. Finally we need an upper bound on
$$e^{k_2s - s/2\b t^2}\E\lcr  \Un_{A_{\a}^r} e^{- s\psi^{}_{\infty}(Y_1^r)/t}\rcr,$$
which will be obtained with the help of large deviation theory. To this aim, we first need to write down properly the stochastic differential equation satisfied by the process $\lacc (Y^r_u, \Phi^r_u), \; u\in[0,1]\racc$, on the event $\lacc T^{\l}_r\ge 1\racc$. Actually, we will have to consider the more general process $\lacc Z^r_u =\lpa Y^r_u, U^r_u, \theta^r_u, \Phi^r_u\rpa, \; u\in [0,1]\racc$, where $U^r$ is the stochastic parallel transform along $Y^r$, and $\theta^r$ is defined by $\theta^r_u = \lpa U^r_u \rpa^*\Theta^r_u$ for every $u\in [0,1]$. 
If $\lacc W_u, \; u\ge 0\racc$ is a Brownian motion in $T_zM$, the very definitions of $Y^r$ and $\Phi^r$ yield first the following equations: 
$$d_\Ito^\n Y^r_u\;=\;\sqrt{r}U^r_u\,dW_u \qquad \hbox{and} \qquad d\Phi_u^Y\;=\;\b t\sqrt{r}(\theta^r)^{\ast} \, dW_u,$$
with the initial conditions $Y^r_0 = z$ and $\Phi^r_0 = 0$. On the other hand, from the definition of the covariant derivative, the equation for $\theta^r$ is given by 
$$ d\theta^r_u\; = \; -\f{r}{2}\lpa U^r\rpa^*\Ric^\sharp(U^r\theta^r)\,du$$
with initial condition $\theta^r_0 = \Id$. It remains to derive the equation for $U^r$. In local coordinates, it is given by
$$\d U^r_u\;=\;-\G(Y^r_u)(\d Y^r_u,U^r_u),$$
where $\d U^r_u$ is the Stratonovich differential of $U^r$, and $\G$ is the
Christoffel symbol of the Laplace-Beltrami connection over $M$ - see e.g. Formula (8.12) in \cite{Emery:89}. The corresponding It\^o equation is 
$$dU^r_u\;=\;-\G(Y^r_u)(dY^r_u,U^r_u)-1/2\G(Y^r_u)(dY^r_u,dU^r_u)-1/2d\G(Y^r_u)(dY^r_u)(dY^r_u,U^r_u)$$
and since in local coordinates,
$$dY^r_u\; =\;\sqrt{r}U^r_u\,dW_u-\f{r}2\trace\G(Y)\,du$$
we finally get 
 \begin{align*}
dU^r_u\;=\;-\sqrt{r}\G(Y^r_u)(U^r_udW_u,U^r_u)&+r/2\Big[ \G(Y^r_u)(U^r_udW_u,\G(Y^r_u)(U^r_udW_u,U^r_u))\\&-d\G(Y^r_u)(U^r_udW_u)(U^r_udW^r_u,U^r_u)\Big].
\end{align*}
Hence, by local boundedness of the Christoffel symbols and their derivatives, we deduce that on $\lacc T^{\l}_r\ge 1\racc$ the process $\lacc Z^r_u, \; u\in[0,1]\racc$ solves in local coordinates an equation of the form 
$$dZ_u^r\; =\;\sqrt{r}\s(Z^r_u)\,dW_u\, +\, rb(Z_u^r)\,du,$$
where $\s$ and $b$ are bounded and uniformly Lipschitz, and whose starting point is $(z,\Id,\Id,0)$. Considering the process $\lacc (r, Z^r_u), \; u\ge 0\racc$ and applying Theorem 5.6.12 in \cite{Dembo-Zeitouni:93}, we see that $\lacc Z^r_u, \; u\in[0,1]\racc$ satisfies a Large Deviation Principle with good rate function
$$I(\phi)\; =\; \inf_{g\in\H^1_0 / \phi = \phi(g)}\lacc \f12\int_0^1\left|\dot g_u\right|^2\,du\racc,$$
where $\H^1_0$ is the usual Cameron-Martin space over $T_zM$ and $\phi(g)$ is the continuous path $\lacc \phi_u =(\phi^1_u,\phi^2_u,\phi^3_u,\phi^4_u), \; u\in [0,1]\racc$ defined followingly: for every $u\ge 0$, $\phi^4_u=\b tg_u$, $\phi^3_u=\Id : T_zM\to T_zM$, and
$$\lacc\begin{array}{l}
\di\phi^1_u\;=\;\int_0^u\phi^2_w\dot g_w\,dw\\
\di\phi^2_u\; =\; \Id-\int_0^u\G(\phi^1_w)(\phi^2_w\dot g_w,\phi^2_w)\,dw
\end{array}\right.$$
The path $u\mapsto \phi_u^1$ is the development on $\{x\in M, \; \rho(z,x) \le \l\}$ in local coordinates of the Euclidean path $u\mapsto g_u$, and $u\mapsto \phi_u^2$ is the parallel transport along $\phi^1$. To obtain an upper bound on 
$$\E\lcr  \Un_{A_{\a}^r} e^{- s\psi^{}_{\infty}(Y_1^r)/t}\rcr,$$
we will apply Theorem 4.3.1 (Varadhan's Lemma) and its consequence Exercise 4.3.11 in \cite{Dembo-Zeitouni:93}. Indeed, the condition (4.3.2) therein on the functional $\psi^{}_{\infty}(Y_1^r)$ is obviously fulfilled, because $\psi^{}_{\infty}(Y_1^r)$ is bounded by a deterministic constant on $\lacc T^{\l}_r\ge 1\racc$. This yields
\begin{align}\label{Varadhan}
\limsup_{s\to \infty}s^{-1}\log \E\lcr  \Un_{A_{\a}^r} e^{- s\psi^{}_{\infty}(Y_1^r)/t}\rcr\; \le \; -\inf_{\phi^4\in A_{\a}}\lpa \b I(\phi^4)+\psi^{}_{\infty}(\phi_1^1)/t\rpa, 
\end{align}
with the above notations and 
$$I(\phi^4)\; =\; 1/2\int_0^1\left|\dot g_u\right|^2\,du.$$
Besides, we deduce from the proof of Lemma 4.3.6 in \cite{Dembo-Zeitouni:93} - which is the main argument to obtain the limsup in Exercise 4.3.11, and hence our above (\ref{Varadhan}) - that the above limsup is uniform in $\alpha$. Indeed, with the notations therein but replacing their $\alpha$ by $\xi$ to avoid confusion, we see that the finite cover of the compact set $\Psi_I(\xi)$ can be chosen independently of our $\alpha$, so that replacing their $\XX$ by our $A_{\a}^r$, the speed of convergence in (\ref{Varadhan}) is dominated by that of the large deviation upper bound for the $\overline{A_{x_i}}'s$ with $x_i\in A_{\a}^r$, hence by that of (\ref{Varadhan}) for $\alpha =0$, which gives the uniformity.

Since $u\mapsto \phi_u^1$ is the development of $u\mapsto g_u$, its Riemannian arclength is given by $\sqrt {2I(\phi^4)}$. But among paths with the same endpoints, geodesics minimize arclength, so that in the above infimum we can consider only paths $\phi$ such
that $g$ is of the form $u\mapsto wu$ with $w\in T_zM$. Now set $A=\sqrt {2I(\phi^4)}$, $e=1/(\b t)$, $a=\|\phi_1^4-v\|$ and $D=\rho(\phi_1^1,
\exp_z(ev))$. With these notations, $\phi_1^4= e^{-1}g_1= e^{-1}w$ and $\phi_1^1=\exp_z w$, so that from the negative curvature of $M$,
$$D\;\ge\;ea.$$ 
We want to find a lower bound for $A^2+2e\psi^{}_{\infty}(\phi_1^1)+e^2$. 
If $\psi^{}_{\infty}(\phi_1^1)\ge 0$, then 
$$A^2+2e\psi^{}_{\infty}(\phi_1^1)+e^2\ge A^2+e^2\ge(A+e)^2/2\ge D^2/2$$
by triangle inequality. If $\psi^{}_{\infty}(\phi_1^1)<0$, then setting $B=\rho(z,\pi(\phi_1^1))$ - recall that $\pi$ is the orthogonal projection on the geodesic $\varphi(x,y)$ - we have  $\psi^{}_{\infty}(\phi_1^1)\ge -B$ by convexity of $\psi_{\infty}$, so that we are left to bound $A^2-2eB+e^2$ from below. Setting $h=\rho\lpa\pi(\phi_1^1),\phi_1^1\rpa$ and using again a comparison theorem together with the negative curvature of $M$, we get $A^2\ge B^2+h^2,$ whence 
$$A^2-2eB+e^2\ge h^2+(B-e)^2\ge \lpa h+|B-e|\rpa^2/2\ge D^2/2$$
by triangle inequality. This entails finally
$$A^2+2e\psi^{}_{\infty}(\phi_1^1)+e^2\ge D^2/2\ge e^2a^2/2$$
and, recalling that $a=\|\phi_1^4-v\|\ge\a$ if $\phi^4\in A_\a$, 
$$1/2\b^2 t^2 + \inf_{g/\phi^4\in A_{\a}}\lpa I(\phi^4)+\psi^{}_{\infty}(\phi_1^1)/\b t\rpa\;\ge\; e^2\a^2/4\;\ge\;\a^2/4r\b^2t^2.$$
Transferring this to the large deviation estimate (\ref{Varadhan}) entails that for every $\e >0$, 
$$\log \E\lcr  \Un_{A_{\a}^r} e^{- s\psi^{}_{\infty}(Y_1^r)/t}\rcr\; \;\le\; s(\e + k_2 - \a^2/4\b t^2)$$
for $s$ big enough, uniformly in $\alpha$. Hence, using (\ref{expect}) and setting $k_2^{\e} = k_2 + \e$, we get 
\begin{eqnarray*}
\E\lcr \Un_{\lacc T^{\l}_r\ge 1\racc} \left\| \Phi^r_1 -v\right\| \f{p_{t- r}(Y^r_1,y)}{p_t(z,y)}\rcr &\le & K\lpa \int_0^\infty 1\wedge
e^{s(k_2^{\e}-  \a^2/4\b t^2)}\,d\a\rpa\\
&\le & K \lpa 2t\sqrt{k_2^{\e}\b} \,+ \,\int_{2t\sqrt{k_2^{\e}\b}}^\infty \!\! e^{s(k_2^{\e}-  \a^2/4\b t^2)}\,d\a\rpa\\
&\le & K \lpa 2\sqrt{k_2^{\e}\b} \,+ \,\int_0^\infty e^{-s\a^2/4\b t^2}\,d\a\rpa\\
&\le & K \lpa 2\sqrt{k_2^{\e}\b} \,+ \,t\sqrt{\pi\b s^{-1}}\rpa,
\end{eqnarray*}
which yields (\ref{finalresult}) and completes the proof of the Theorem. 

\fin

Unfortunately, the above theorem is not sufficient to entail (\ref{E:9}), because there we strongly need uniformity on $\omega \in E(s)$ and $t\in[0,1/2]$. On the other hand, it seems difficult to provide a uniform version of Theorem \ref{heat3} without further assumption on the curvature tensor of $M$: in order to apply Theorem 5.6.12 and Exercise 4.3.11 in \cite{Dembo-Zeitouni:93}, it is necessary to have {\em uniform} boundedness and Lipschitz properties for $\sigma$ and $b$, and this uniformity fails whenever the Christoffel symbols or their derivatives up to order two are not bounded on a fixed neighbourhood of the origin of the exponential maps. This situation is possible, as can be seen from the example of a two-dimensional radially symmetric manifold with prescribed sectional curvature $\kappa(r) = -(1 +\cos^2r^2)$, where the second derivative of the Christoffel symbols at the origins of the  
exponential maps is not bounded (we leave the details to the reader).

For this reason we need Assumption \ref{Christo} on M to obtain the following uniform version of Theorem \ref{heat3} which, combined with Lemma \ref{Laplace2}, Proposition \ref{E2(s)}, and the end of the proof in Section \ref{Rankone}, will finish the proof of Theorem \ref{T1}.

\begin{thm} \label{heat4} Under Assumption \ref{Christo} and with the notations of Theorem \ref{heat3}, for every $\e \in ]0,1]$
$$\rho(z,y)^{-1}\grad\log p_t(\newdot,y)(z)\; \to \; t^{-1}\dot\varphi(z,y)(0)$$
when $\rho(z,y)\to +\infty$, uniformly on $t\in[\e,1]$, $z\in M$, and $v\in T_zM$.
\end{thm}

\begin{proof} The proof follows almost verbatim from that of Theorem \ref{heat3}, since Lemma \ref{cLambda} is already uniform in $z$. All we need to check is that the limsup in (\ref{Varadhan}) is uniform in $\alpha > 0$ {\em and} $z\in M$. First, Assumption \ref{Christo} yields the desired uniform bound for $b$ and $\sigma$ and their derivatives, and then we can check from the proof of Theorem 5.6.12 in \cite{Dembo-Zeitouni:93} that the speed of convergence in the Large Deviation upper bound for $\lacc Z^r_u, \; u\in[0,1]\racc$ does not depend on $z$. The same holds concerning Lemma 4.3.11 in \cite{Dembo-Zeitouni:93}, but the argument is a bit more subtle because the compact set $\Psi_I(\xi)$ in Lemma 4.3.6  in \cite{Dembo-Zeitouni:93} - with our notation for $\xi$ instead of their $\alpha$ to avoid confusion - also depends on $z$. However, the uniformity in $z$ remains as far as the large deviation upper bound for the $\overline{A_{x_i}}'s$ is concerned, because of the uniform LDP for $\lacc Z^r_u, \; u\in[0,1]\racc$. This completes the proof.

\end{proof}

\begin{remark} During the proofs of Theorems \ref{heat3} and \ref{heat4}, another method for majorizing
$$E\; :=\;\esp\left[\Un_{\{T_r^\l\ge 1\}}\left\|\Phi_1^r-v\right\|\f{p_{t-r}(Y_1^r,y)}{p_t(z,y)}\right]$$
could be possible, not relying on Varadhan's Lemma. We first decompose
 \begin{align*}
E\; =\;& \esp\left[\Un_{\{T_r^\l\ge
    1\}}\Un_{\left\{\left\|\Phi_1^r-v\right\|\le
    2t\sqrt{k_2\beta}\right\}}\left\|\Phi_1^r-v\right\|\f{p_{t-r}(Y_1^r,y)}{p_t(z,y)}\right]\\
&\qquad\qquad+\esp\left[\Un_{\{T_r^\l\ge
    1\}}\Un_{\left\{\left\|\Phi_1^r-v\right\|>
    2t\sqrt{k_2\beta}\right\}}\left\|\Phi_1^r-v\right\|\f{p_{t-r}(Y_1^r,y)}{p_t(z,y)}\right]\\
\le\; &2t\sqrt{k_2\beta}\; +\; \esp\left[\Un_{\{T_r^\l\ge
    1\}}\Un_{\left\{\left\|\Phi_1^r-v\right\|>
    2t\sqrt{k_2\beta}\right\}}\left\|\Phi_1^r-v\right\|\f{p_{t-r}(Y_1^r,y)}{p_t(z,y)}\right],
\end{align*}
so that we need an upper bound on
$$E'\; :=\;  \esp\left[\Un_{\{T_r^\l\ge 1\}}\Un_{\left\{\left\|\Phi_1^r-v\right\|> 2t\sqrt{k_2\beta}\right\}}\left\|\Phi_1^r-v\right\|\f{p_{t-r}(Y_1^r,y)}{p_t(z,y)}\right].$$
Since
$$\f{p_{t-r}(Y_1^r,y)}{p_t(z,y)} \;\le\; Ke^{k_2s-\f{s}{2\b t^2}-\f{s}{t}\psi_\infty(Y_1^r)}$$
on the event $\{T_r^\l\ge 1\}$, we first have 
\begin{equation}
\label{E}
E\;\le\; Ke^{k_2s}e^{-\f{s}{2\b t^2}}\esp\left[\Un_{\{T_r^\l\ge 1\}}\Un_{\left\{\left\|\Phi_1^r-v\right\|> 2t\sqrt{k_2\beta}\right\}}\left\|\Phi_1^r-v\right\|e^{-\f{s}{t}\psi_\infty(Y_1^r)}\right].
\end{equation}
Let $\lacc Z^r_t,\; t\ge 0\racc$ be the solution to the SDE
$$Z^r_u\; =\; Y^r_u\;  -\; \f{1}{\b t}\int_0^u\grad\psi_\infty(Z_w^r)\,dw,$$
where $\lacc Y^r_t,\; t\ge 0\racc$ is as above a Brownian motion starting from $z$ and with speed $r$. Considering the martingale
$$N_u\;=\;\f{1}{r\b t}\int_0^u\left\langle \grad\psi_\infty(Z_w^r),d_\Ito^\n
Z_w^r+\f{1}{\b t}\grad\psi_\infty(Z_w^r)\,dw\right\rangle,$$
we see that under $\QQ=\exsp(N)\cdot \PP$, $Z^r$ has the same law as $Y^r$. Hence, the right-hand side of \eqref{E} can be rewritten as 
$$Ke^{k_2s}e^{-\f{s}{2\b t^2}}\esp\left[\Un_{\left\{T_r^\l(Z^r)\ge 1,\; \left\|\Phi_1^r(Z^r)-v\right\|> 2t\sqrt{k_2\beta}\right\}}\left\|\Phi_1^r(Z^r)-v\right\|e^{-\f{s}{t}\psi_\infty(Z_1^r)}\exsp(N)_1\right]$$
where $T_r^\l(Z^r)$ and $\Phi_1^r(Z^r)$ are defined in replacing $Y^r$ by $Z^r$. Noticing that $\|\grad\psi_\infty\|=1$, we have
$$\exsp(N)_1=e^{\f{1}{2r\b^2t^2}}\exp\left(\f{1}{r\b t}\int_0^1\left\langle\grad\psi_\infty(Z_w^r),d_\Ito^\n Z_w^r\right\rangle\right).$$
Applying It\^o's formula yields
$$\psi_\infty(Z_1^r)=\int_0^1\left\langle \grad\psi_\infty(Z_u^r),d_\Ito^\n
Z_u^r\right\rangle+\f{r}2\int_0^1\Delta \psi_\infty(Z_u^r)\,du,$$
so that finally 
\begin{eqnarray*}
E & \le & Ke^{k_2s} \esp\left[\Un_{\left\{T_r^\l(Z^r)\ge 1,\; \left\|\Phi_1^r(Z^r)-v\right\|> 2t\sqrt{k_2\beta}\right\}}\left\|\Phi_1^r(Z^r)-v\right\|e^{-\f{1}{\b t}\int_0^1\Delta\psi_\infty(Z_u^r)\,du}\right]\\
& \le & Ke^{k_2s} \esp\left[\Un_{\left\{T_r^\l(Z^r)\ge 1,\; \left\|\Phi_1^r(Z^r)-v\right\|> 2t\sqrt{k_2\beta}\right\}}\left\|\Phi_1^r(Z^r)-v\right\|\right]\\
& \le & Ke^{k_2s} C(\b)\sqrt{r} \,\pb\lcr T_r^\l(Z^r)\ge 1,\; \left\|\Phi_1^r(Z^r)-v\right\|> 2t\sqrt{k_2\beta}\rcr^{1/2},
\end{eqnarray*}
where in the second line we used the fact (coming from the convexity of $\psi_\infty$) that $\Delta \psi_\infty \ge 0$ and where in the third line we have majorized the $L^2$-norm of $\left\|\Phi_1^r(Z^r)-v\right\|$. Using the same Large Deviation Principle for $Z^r$ as above, one can prove that
$$\pb\lcr T_r^\l(Z^r)\ge 1,\; \left\|\Phi_1^r(Z^r)-v\right\|> 2t\sqrt{k_2\beta}\rcr\; \le \; e^{-k_{\beta}s}$$
as $s\to +\infty$, with
$$k_{\beta} \; = \; \inf_{\tiny \begin{array}{c} g\in\H^1_0, \phi = \phi(g)\\\lva\lva \phi^4_1 - v\rva\rva \ge 2t\sqrt{k_2\beta}\end{array}}\lacc \f12\int_0^1\lva\lva\phi^2_u\dot g_u + \frac{1}{\beta t} \grad\psi_{\infty}(\phi^1_u)\rva\rva^2\,du\racc.$$
However, because of the difficult tractability of $\grad\psi_{\infty}$, except in the trivial flat case we got stuck in proving that $k_{\beta}\to +\infty$ as $\beta\to 0,$ which would be  enough to complete the proof of Theorems \ref{heat3} and \ref{heat4}. The advantage of this method is that there would be no Busemann function anymore under the integral for the large deviation estimates, so that we would only need to use Theorem 5.6.12 in \cite{Dembo-Zeitouni:93}. 

\end{remark}

%
%

\noindent
{\em Acknowledgement:} Part of this work was done during a visit at the Technische Universit\"at Berlin of the second-named author, who would like to thank Michael Scheutzow for his kind hospitality and the excellent working conditions.

%
%

\providecommand{\bysame}{\leavevmode\hbox to3em{\hrulefill}\thinspace}


\begin{thebibliography}{10}

\bibitem{Anker-Ostellari:03}
J.~P.~Anker and P.~Ostellari, \emph{The heat kernel on non compact symmetric spaces.} In: \emph{Lie groups and symmetric spaces}, Amer. Math. Soc. Transl. Ser. 2, \textbf{210}, Amer. Math. Society, Providence, RI (2003), 27--46.

\bibitem{Berger:03}
M.~Berger, \emph{A Panoramic View of Riemannian Geometry.} Springer-Verlag, 2003.

\bibitem{Bismut:84}
J.~M.~Bismut, \emph{Large Deviations and the Malliavin Calculus.} Birkh\"auser, 1984.

\bibitem{Bridson-Haefliger:99}
M.~R.~Bridson and A.~Haefliger, \emph{Metric Spaces of Non-Positive Curvature}, Springer-Verlag, Berlin, 1999. 

\bibitem{Cheeger-Ebin:75}
J.~Cheeger and D.~G.~Ebin, \emph{Comparison Theorems in Riemannian Geometry}, North-Holland, 1975.

\bibitem{Dembo-Zeitouni:93}
A.~Dembo and O.~Zeitouni, \emph{Large Deviations Techniques and Applications}, Jones and Barlett Publishers, Boston, 1993. 

\bibitem{Eberle:02}
A.~Eberle, \emph{Absence of spectral gaps on a class of loop spaces}, J. Math. Pures Appl. (9) \textbf{81} (2002), 915--955.

\bibitem{Emery:89}
M.~Emery, \emph{Stochastic Calculus in Manifolds}, Springer-Verlag, Berlin, 1989.

\bibitem{Feller:54}
W.~Feller, \emph{Diffusion processes in one dimension}, Trans. Amer. Math. Soc. \textbf{77} (1954), 1--31.

\bibitem{Carlo-Maria:93}
S.~Giulini and W.~Woess, \emph{The Martin compactification of the Cartesian product of two hyperbolic spaces}, J. Reine Angew. Math. \textbf{444}  (1993), 17--28.

\bibitem{Helgason:78}
S.~Helgason, \emph{Differential Geometry, Lie Groups and Symmetric Spaces}, Academic Press, 1978. 

\bibitem{Helgason:84}
S.~Helgason, \emph{Groups and Geometric Analysis}, Academic Press, 1984. 
 
\bibitem{Hsu:02}
E.~P.~Hsu, \emph{Stochastic Analysis on Manifolds}, Amer. Math. Society, Providence, RI, 2002.

\bibitem{Koornwinder:84}
T.~H.~Koornwinder, \emph{Jacobi functions and analysis on non compact semisimple Lie groups.} In: \emph{Special functions: group theoretical aspects and applications}, R.~A.~Askey et al. (eds.), Reidel (1984), 1--85.

\bibitem{Lohoue-Rychener:82}
N.~Lohoue and T.~Rychener, \emph{Die Resolvente von $\Delta$ auf symmetrischen Ra\"umen von nichtkompakten Typ}, Comment. Math. Helvet. \textbf{57}  (1982),  445--468.

\bibitem{Lorang-Roynette:96}
G.~Lorang and B.~Roynette, \emph{Etude d'une fonctionnelle li\'ee au pont de Bessel}, Ann. Inst. H. Poincar\'e Probab.  Statist. \textbf{32}  (1996),  no.~1, 107--133.

\bibitem{Magnus:66}
W.~Magnus, F.~Oberhettinger and R.~P.~Soni, \emph{Formulas and Theorems for the Special Functions of Mathematical Physics,} Springer-Verlag, New York, 1966.

\bibitem{Mandl:68}
P.~Mandl, \emph{Analytical Treatment of One-Dimensional Markov Processes}, Academia, Prague, and Springer-Verlag, New-York, 1968.

\bibitem{Norris:93}
J.~R.~Norris, \emph{Path integral formulae for heat kernels and their derivatives}, Probab. Theory Related Fields {\bf 94}, 525--541, 1993.
  
\bibitem{Revuz-Yor:99}
D.~Revuz and M.~Yor, \emph{Continuous Martingales and Brownian Motion}, 
 Springer-Verlag, Berlin, 1999. 

\bibitem{Simon:02}
T.~Simon, \emph{Concentration of the Brownian bridge on the hyperbolic plane}, Ann. Probab.  \textbf{30}  (2002),  no.~4, 1977--1989.

\bibitem{Thalmaier:97}
A.~Thalmaier, \emph{On the Differentiation of Heat Semigroups and Poisson Integrals}, Stoch. Stoch. Rep.  \textbf{61} (1997), 297--321.

\bibitem{Thalmaier-Wang:98}
A.~Thalmaier and F.Y.~Wang, \emph{Gradient estimates for harmonic functions on regular domains in Riemannian manifolds}, J. Funct. Anal. \textbf{155} (1998), 109-124.

\end{thebibliography}
\end{document}